\documentclass[11pt,leqno]{article}
\usepackage{ graphicx, amsfonts, amsthm, amsxtra, amsmath,verbatim, multicol}
\usepackage[mathscr]{euscript}
\makeindex
\begin{document}

\begin{center}
{\LARGE\bf Automorphic forms for triangle groups
}
\\
\vspace{.25in} {\large {\sc Charles F. Doran, Terry Gannon,\footnote{
Department of Mathematical and Statistical Sciences,
University of Alberta, Edmonton AB T6G 2G1. } Hossein Movasati and Khosro M. Shokri }}\footnote{
Instituto de Matem\'atica Pura e Aplicada, IMPA,
Estrada Dona Castorina, 110,
22460-320, Rio de Janeiro, RJ, Brazil,}
\end{center}

\newtheorem{theo}{Theorem}
\newtheorem{exam}{Example}
\newtheorem{coro}{Corollary}
\newtheorem{defi}{Definition}
\newtheorem{prob}{Problem}
\newtheorem{lemm}{Lemma}
\newtheorem{prop}{Proposition}
\newtheorem{rem}{Remark}
\newtheorem{conj}{Conjecture}
\newtheorem{exe}{Exercise}
\newcommand\diff[1]{\frac{d #1}{dz}} 
\def\End{{\rm End}}              
\def\hol{{\rm Hol}}
\def\sing{{\rm Sing}}            
\def\spec{{\rm Spec}}            
\def\cha{{\rm char}}             
\def\Gal{{\rm Gal}}              
\def\jacob{{\rm jacob}}          
\def\Z{\mathbb{Z}}                   
\def\O{{\cal O}}                       

\def\C{\mathbb{C}}                   
\def\as{\mathbb{U}}                  
\def\ring{{\sf R}}                             
\def\R{\mathbb{R}}                   
\def\N{\mathbb{N}}                   
\def\A{\mathbb{A}}                   

\def\D{\mathbb{D}}                   
\def\uhp{{\mathbb H}}                
\def\bbH{{\mathbb H}}
\newcommand\ep[1]{e^{\frac{2\pi i}{#1}}}
\newcommand\HH[2]{H^{#2}(#1)}        
\def\Mat{{\rm Mat}}              
\newcommand{\mat}[4]{
     \begin{pmatrix}
            #1 & #2 \\
            #3 & #4
       \end{pmatrix}
    }                                
\newcommand{\matt}[2]{
     \begin{pmatrix}                 
            #1   \\
            #2
       \end{pmatrix}
    }
\def\ker{{\rm ker}}              
\def\cl{{\rm cl}}                
\def\dR{{\rm dR}}                

\def\hc{{\mathsf H}}                 
\def\Hb{{\cal H}}                    
\def\GL{{\rm GL}}                
\def\pedo{{\cal P}}                  
\def\PP{\tilde{\cal P}}              
\def\cm {{\cal C}}                   
\def\K{{\mathbb K}}                  
\def\k{{\mathsf k}}                  
\def\F{{\cal F}}                     
\def\M{{\cal M}}
\def\RR{{\cal R}}
\newcommand\Hi[1]{\mathbb{P}^{#1}_\infty}
\def\pt{\mathbb{C}[t]}               
\def\W{{\cal W}}                     
\def\Af{{\cal A}}                    
\def\gr{{\rm Gr}}                
\def\Im{{\rm Im}}                
\newcommand\SL[2]{{\rm SL}(#1, #2)}    
\newcommand\PSL[2]{{\rm PSL}(#1, #2)}  
\def\Res{{\rm Res}}              

\def\L{{\cal L}}                     
\def\Aut{{\rm Aut}}              
\def\any{R}                          
\newcommand\ovl[1]{\overline{#1}}    

\def\per{{\sf  pm}}  
\def\T{{\cal T }}                    
\def\tr{{\sf tr}}                 
\newcommand\mf[2]{{M}^{#1}_{#2}}     
\newcommand\bn[2]{\binom{#1}{#2}}    
\def\ja{{\rm j}}                 
\def\Sc{\mathsf{S}}                  
\newcommand\es[1]{g_{#1}}            
\newcommand\V{{\mathsf V}}           
\newcommand\Ss{{\cal O}}             
\def\rank{{\rm rank}}                
\def\diag{{\rm diag}}
\def\BM{{\sf H}}
\def\Fi{{S}}
\def\Ra{{\rm R}}

\def\Q{{\mathbb Q}}                   

\def\P{\mathbb{P}}

\def\bG{\mathbf{G}}
\def\bbZ{\mathbb{Z}} 
\def\uhp{\mathbb{H}}
\def\bbR{\mathbb{R}}
\def\bbQ{\mathbb{Q}}
\def\bbC{\mathbb{C}}
\def\bbP{\mathbb{P}}
\def\bbW{\mathbb{W}}
\def\bbX{\mathbb{X}}
\def\I{\mathrm{i}}
\def\mft{\mathfrak{t}}
\def\cP{\mathcal{P}}
\def\cM{\mathcal{M}}
\def\qh{\hat{q}}

\begin{abstract}
For triangle groups, the (quasi-)automorphic forms are known just as explicitly as for the modular
group SL$(2,\bbZ)$. We collect these expressions here, and then interpret them using 
 the Halphen differential equation. We study the arithmetic properties of 
 their Fourier coefficients  at cusps and Taylor coefficients at elliptic fixed-points ---
 in both cases integrality is related to the arithmeticity of the triangle group.  
As an application of our formulas, we provide an explicit modular interpretation of 
periods of 14 families of Calabi-Yau threefolds over the thrice-punctured sphere.
\end{abstract}

{\tiny
\tableofcontents
}
\section{Introduction}

Although modular forms for congruence subgroups of the modular group PSL$(2,\bbZ)=\Gamma(1)$
go back to Euler, modular forms for more general Fuchsian groups (usually called
\textit{automorphic forms}) go back to Poincar\'{e}. He proved their existence
by constructing 
functions (Fuchsian-theta series in his terminology) which nowadays are known as Poincar\'{e} series.
Independently of Poincar\'e, G. Halphen in \cite{ha81, ha81-0} introduced a differential equation in three variables and three parameters, 
which nowadays bears his name. His motivation was a particular case
studied by Darboux in \cite{da78} and he proved that in such a case the differential equation is 
satisfied by the logarithmic derivatives of theta functions. 
Despite the fact that Poincar\'e and Halphen 
were contemporaries and compatriots, the main relation between these works was not clearly understood, and Halphen's contribution was largely forgotten, only to be rediscovered several
times.

The modular forms and functions for the modular group $\Gamma(1)$ have of course been
well understood for many decades. What is less well known is that there is a natural infinite
class of Fuchsian groups --- the so-called \textit{triangle groups} --- where the automorphic forms
and functions can be determined just as explicitly, even though all but a few are incommensurable
with $\Gamma(1)$. 

Let $\Gamma\le \mathrm{PSL}(2,\bbR)$ be any genus-0 finitely generated Fuchsian group of the first kind.
(See the following section for the definitions of these and other technical terms.)
This means that $\Gamma\backslash\uhp_\Gamma$ is topologically a sphere, where $\uhp_\Gamma$ denotes the
upper half-plane $\uhp$ extended by the cusps of $\Gamma$ (if any).
 Let $n_{cp}$ be the number of cusps and $n_{el}$ be the number of
elliptic fixed-points, and write $2\le n_i\le\infty$ for the orders of their stabilizers. Then 
Gauss-Bonnet implies $2<\sum_{j=1}^{n_{cp}+n_{el}}(1-1/n_j)$ (see e.g.
Theorem 2.4.3 of \cite{Miy} for a generalization) and hence we have the
inequality $n_{cp}+n_{el}\ge 3$.  The field of automorphic functions of $\Gamma$ is
$\bbC (J_\Gamma)$ where the generator $J_\Gamma$ maps $\Gamma\backslash\uhp_\Gamma$ bijectively
onto the Riemann sphere $\bbP^1$. Knowing such a uniformizer $J_\Gamma$ determines
explicitly (in principle) all automorphic and quasi-automorphic forms. If $\Gamma$ is commensurable with
$\Gamma(1)$ (i.e. when $\Gamma\cap\Gamma(1)$ has finite index in both
$\Gamma$ and $\Gamma(1)$), then (in principle) a generator $J_\Gamma$ can be determined from
e.g. the Hauptmodul $j(\tau)=q^{-1}+196884q+\cdots$ of $\Gamma(1)$. When $\Gamma$
is not necessarily commensurable, it is useful to recall that $J_\Gamma$ will satisfy a nonlinear
third order differential equation
\begin{equation}
-2{J_{\Gamma}'''(\tau)\over J_{\Gamma}'(\tau)}+{3}{J_{\Gamma}''(\tau)^2\over J_{\Gamma}'(\tau)^2}
=J_{\Gamma}'(\tau)^2Q_\Gamma(J_\Gamma(\tau))\label{schwder}\end{equation}
coming from the Schwarzian derivative, where the prime here denotes ${{\rm d}\over{\rm d}\tau}$. 

The Schwarzian equation \eqref{schwder} is rather complicated. It can be replaced
by a much simpler system of first order differential equations in $n_{cp}+n_{el}$ variables,
subject to $n_{cp}+n_{el}-3$ quadratic (nondifferential) constraints. In this generality,
the result is due to Ohyama \cite{Oh}, but the key ideas go back to the 19th century.
In particular, Halphen \cite{ha81} associated the system
\begin{equation}
\label{tatlis}
\left \{ \begin{array}{l}
 t_1'=
( a-1)(t_1t_2+t_1t_3-t_2t_3)+( b+ c-1)t_1^2\\
t_2' =
(b-1)(t_2t_1+t_2t_3-t_1t_3)+( a+ c-1)t_2^2\\
 t_3' =
( c-1)(t_3t_1+t_3t_2-t_1t_2)+(a+ b-1)t_3^2
\end{array} \right.  \ \ ,
\end{equation}
where the prime denotes $d/d \tau$, to Gauss' hypergeometric equation
\begin{equation}
\label{gauss}
z(1-z)y''+(a+c-(a+b+2c)z)y'-(a+b+c-1)cy=0\,, 
\end{equation}
where now the prime denotes $d/{dz}$, and Brioschi \cite{Bri} showed its equivalence to the corresponding version of \eqref{schwder}
(namely, \eqref{schw} below). The Halphen system
\eqref{tatlis} has been rediscovered several times (including by one of the authors of this
 paper!), and over the past century has
appeared in the study of monopoles, self-dual Einstein equations, WDVV equations, 
mirror maps, etc. 
In \cite{HarnMcK} the authors have used solutions of Halphen equation for many particular cases, including those with an arithmetic
triangle group, to obtain replicable uniformizations of punctured Riemann surfaces of genus zero.  Further particular cases of Halphen
equation solved by  classical theta series or modular forms are discussed in \cite{ab06}. The idea to use Halphen equation and find
new automorphic forms seems to be neglected in the literature.

Now, $Q_\Gamma(z)$ in \eqref{schwder} is a rational function depending on $n_{cp}+n_{el}-2$ parameters.
Unfortunately, these parameters depend on $\Gamma$ in a very complicated nonalgebraic
way and in general closed formulae for them cannot be found (see e.g. \cite{ven} for
an analysis of this question). 
However, when $n_{cp}+n_{el}=3$ (the minimum value possible), this single parameter
can be determined explicitly, using classical results on hypergeometric functions.
In this case --- where $\Gamma$ is a {triangle group} --- $J_\Gamma(\tau)$
and hence all quasi-automorphic forms for $\Gamma$ can be explicitly determined.

One of the purposes of this paper is to write these explicit expressions down. 
Special cases and partial results are scattered throughout the literature, but
to our knowledge these expressions haven't appeared in the literature with this explicitness
and in this generality, and certainly not all in one place.

We do this in two ways. We begin with the classical approach, because of its familiarity: the 
 multivalued ratio $\tau(z)$ of
two solutions to the hypergeometric equation can in certain circumstances be regarded
as the functional inverse of an automorphic function $z(\tau)$ for a triangle group. This 
determines $z(\tau)$ completely, but it is convenient to use \eqref{schwder} to recover its
$q$-expansion. Differentiating $z(\tau)$ once yields all automorphic forms; differentiating
it a second time yields all quasi-automorphic forms. Although
the basic ideas of this derivation are classical, going back to Fuchs and Poincar\'e, the details are unpleasant.
\textit{Our second approach,  using the
Halphen equation, is independent and turns this on its head,} even though the underlying
mathematics is again that of the hypergeometric equation. We interpret solutions of
Halphen's equation, when lifted to $\uhp$,  as quasi-automorphic  forms for a triangle
group. Taking differences yields all automorphic forms, and ratios then yield all automorphic
functions.  

We suggest that in most respects, the (quasi-)automorphic forms of the triangle groups 
are close cousins of those of the modular group and can be studied analogously,
even though these groups are (usually) not commensurable with $\Gamma(1)$ (and so e.g.
Hecke operators cannot be applied). In particular,
everything is as explicit for arbitrary triangle groups as it is for the modular group.

Now, when the group contains a congruence subgroup $\Gamma(n)$ of $\Gamma(1)$,
such modular forms have many arithmetic  properties. It is natural to ask whether any such
arithmeticity survives for general triangle groups. We explore the arithmeticity of both the
local expansions at cusps and at elliptic fixed-points. The latter expansions are far less familiar, even
though they were familiar to e.g. Petersson in the 1930s \cite{Pet}, but they deserve 
more attention than they have received.  For example,
 Rodriguez Villegas--Zagier \cite{RVZ} interpret the expansion coefficients of
the Dedekind eta $\eta(\tau)$ at $\omega=e^{2\pi\I/3}$ in terms of central values of Hecke
L-functions.

The triangle groups are extremely special among the Fuchsian groups for a number of reasons,
for instance:

\begin{itemize}
\item[(i)]
One is a consequence of Belyi's theorem. A Fuchsian group is a subgroup of finite index in a 
triangle group, iff for each weight $k\in2\bbZ$, there is a basis of the $\bbC$-space of weight-$k$
holomorphic automorphic forms whose expansion coefficients are all algebraic numbers (see
e.g. \cite{scholl}).
Of course, these coefficients are the primary reason for the importance of any automorphic forms.

\item[(ii)] The complement of a knot in $S^3$ has universal cover $\widetilde{\mathrm{SL}}(2,\bbR)$
(the universal cover of SL$(2,\bbR)$), iff the knot is a torus knot \cite{RV}. In particular, the
$(p,q)$-torus knot is diffeomorphic to $\widetilde{\mathrm{SL}}(2,\bbR)/G$ for a certain
lift of the $(p,q,\infty)$-triangle group. For example, the complement of the trefoil is 
  $\widetilde{\mathrm{SL}}(2,\bbR)/\widetilde{\mathrm{SL}}(2,\bbZ)$.
 The relevance to this here is  that an automorphic form, of arbitrary weight, for $\Gamma$ lifts to a
function on $\widetilde{\mathrm{SL}}(2,\bbR)/\widetilde{\Gamma}$. The relevance to
torus knots of the automorphic forms of the $(p,q,\infty)$-triangle group is developed in
\cite{Tsa}, following \cite{Mil} and Section 2.4.3 of \cite{Ganb}. Now, recall that Gopakumar--Vafa duality would imply that the Chern-Simons
knot invariants arise as Gromov--Witten invariants.  This has been verified explicitly in \cite{BEM}
for
the torus knots, by independently computing the two sets of invariants and showing they are
equal. It seems very possible that reinterpreting \cite{BEM} using automorphic forms for
triangle groups would at least simplify their calculation, and could lead to a more conceptual
explanation of the equality. 

\item[(iii)] We see below that periods of some Calabi-Yau three-folds with 1-dimensional
moduli spaces can be interpreted as vector-valued automorphic forms for certain 
triangle groups (e.g. $(5,\infty,\infty)$ for the dual of the quintic). Independently,
 all 26 sporadic finite simple groups are quotients of certain triangle groups \cite{wil},
e.g. the Monster is a quotient of ${(2,3,7)}$ (and hence $\Gamma(1)$).
This implies that, for each sporadic group $G$, there will exist vector-valued automorphic
forms  for
some triangle group, whose multiplier $\rho$ factors through to a faithful representation of $G$.
\end{itemize}

In \cite{ho07-1}, the author (HM) derived the Halphen differential equation using the inverse of a period map. 
One advantage of this point of view is the introduction of modular-type forms for finitely generated subgroups of PSL$(2,\C)$ which may 
not be even discrete, something which must sound dubious to most number theorists. 
Since \cite{ho07-1} focusses on the differential  and geometric aspects of such modular-type forms, we felt that 
we should now look at number theoretic aspects. The triangle groups provide interesting but 
nontrivial toy models,
 where the group is discrete but the automorphic forms are not 
so well-studied. This text is partly a result of this effort.
We find it remarkable how naturally the (quasi-)automorphic
forms for triangle groups arise in the Halphen system \eqref{tatlis}. We believe this
observation is new (at least in this generality). In this case, the parameters $a,b,c$
must be rational --- in fact the combinations $1-a-b,1-c-b,1-a-c$ will equal the angular parameters $v_i=1/m_i$, for $i=1,2,3$ respectively, where $m_i\in\bbZ_{>0}\cup\{\infty\}$. 
However, 
 some sort of modularity appears to persist though even when these angular parameters are complex.

Our main motivation for writing  this paper is to establish the background needed to
understand the modularity of the mirror map
for examples such as the Calabi-Yau quintic, by relating the Halphen approach of one
of the authors with that of vector-valued automorphic forms of another author.
This required having completely explicit descriptions of the automorphic forms for the triangle group $(5,\infty,\infty)$, and as we couldn't find this adequately treated  in the literature we 
did the calculations ourselves. The application to mirror maps will be forthcoming,
although an initial step is provided in Section 6.

The outline of the paper is as follows. 
Section 2 provides the
classical (i.e. hypergeometric) calculation of all data for the automorphic forms of the triangle 
groups.  
Section 3 recovers this data using solutions to Halphen's equation; we believe this approach
is new.  Section 4 specializes to 
the triangle groups commensurable with the modular group.
Section 5 explores the arithmeticity of the Fourier and Taylor coefficients.  Section 6 applies this
material to periods of Calabi-Yau threefolds. Our proofs are collected in Section 7. 
Relevant facts on hypergeometric functions are collected in Appendix A.

Here is some notations used throughout the text.
\begin{itemize}
\item
$\mft=(m_1,m_2,m_3)$: triangle group type;
\item
$\uhp$ resp. $\uhp_\mft$: the upper half-plane resp. extended upper half-plane;
\item
$\Gamma_{\mft}\subset \SL 2\R$: the realization of the triangle group of type $\mft$;
\item
$\gamma_i$, $i=1,2,3$: matrix generators of $\Gamma_\mft$ (see \eqref{trigen});\item
$\zeta_i$, $i=1,2,3$: fixed-points of  $\gamma_i$ (see \eqref{corners});\item
$q_i$ resp. $\tilde{q}_i$: the local coordinate resp. normalized local coordinate, at $\zeta_i$;
\item
$J_{\mft}$: the normalized Hauptmodul associated to the group $\Gamma_\mft$ (see \eqref{Jnorm});\item
$v_i=\frac{1}{m_i},\ i=1,2,3$: the angular parameters; \item
$(a,b,c)$ resp. $(\tilde{a},\tilde{b},\tilde{c})$:  parameters of the Halphen resp. hypergeometric systems;\item
$(t_1,t_2,t_3)$: the solution of the Halphen system, defined in \S3.
\end{itemize}

\section{Classical computation of (quasi-)automorphic forms}
\label{24oct02}
In this section we give the classical approach for computing automorphic forms through the Schwarzian and hypergeometric differential equations.

\subsection{Background}

See e.g.
\cite{Miy} for the basics of Fuchsian groups and their automorphic forms. A \textit{Fuchsian group} $\Gamma$ is a discrete subgroup of PSL$(2,\bbR)=\mathrm{SL}(
2,\bbR)/\{\pm 1\}$, the group of orientation-preserving isometries of the upper half-plane
$\bbH:=\{x+iy\mid y>0\}$. $\Gamma$ is called \textit{of first class} (the class of primary interest) if its fundamental
domains in $\bbH$ have finite hyperbolic area. $\gamma\in\Gamma$ is called \textit{parabolic}
if $\gamma$ has precisely one fixed-point on the boundary $\partial\bbH=\bbR\bbP^1=
\bbR\cup\{\I\infty\}$;
$x\in\bbR\cup\{\I\infty\}$ is called a \textit{cusp} of $\Gamma$ if it is fixed by some parabolic
$\gamma\in\Gamma$. The \textit{extended half-plane}
together with all cusps; then for $\Gamma$ of first class, the orbits $\Gamma\backslash\bbH_\Gamma$
naturally form a compact surface. The genus of this surface is called the genus of $\Gamma$.

If $\I\infty$ is a cusp of $\Gamma$, we call the smallest $h>0$ with $\gamma_{\infty;h}:=\left({1\atop 0}{h\atop 1}\right)
\in\Gamma$ the \textit{cusp-width} $h_\infty$. If $x\in\bbR$ is a cusp, its cusp-width $h_x$
is the smallest $h>0$ for which  $\gamma_{x;h}:=\left({0\atop 1}{-1\atop -x}\right)^{-1}\left({1\atop 0}{h\atop 1}\right)
\left({0\atop 1}{-1\atop -x}\right)\in\Gamma$.
The other special points in $\bbH_\Gamma$ are the \textit{elliptic fixed-points}, which are $z\in\bbH$
stabilized by a nontrivial $\gamma\in\Gamma$. For each $z=x+\I y\in\bbH$, the stabilizer in $\Gamma$
is finite cyclic, generated by
$$\gamma_{z;n}:=\left({y^{-1/2}\atop 0}\,{-y^{-1/2}x\atop y^{1/2}}\right)^{-1}\left({\cos(\pi/n)\atop 
-\sin(\pi/n)}{\sin(\pi/n)\atop \cos(\pi/n)}\right)
\left({y^{-1/2}\atop 0}\,{-y^{-1/2}x\atop y^{1/2}}\right)$$
for a unique positive integer $n=n_z$ called the \textit{order} of $z$.  Write $n_x=\infty$ for a 
cusp $x$.

These numbers $h_x,n_z$ are clearly constant along $\Gamma$-orbits. Let $n_{el}$ denote the
number of $\Gamma$-orbits of elliptic fixed-points, and  $n_{cp}$
 the number of $\Gamma$-orbits of cusps. Both $n_{el}$ and $n_{cp}$ must be finite,
but can be zero; moreover, $n_{el}+n_{cp}\ge 3$. 

For $z\in \uhp_\Gamma$, define M\"obius transformations $\tau\mapsto \tau_{z}$, local coordinates 
$q_{z}$ and  automorphy factors
$j_{z}(k;\tau)$ as follows. Choose $\tau_{\infty}=\tau$, $q_{\infty}
=e^{2\pi\I\tau/h_\infty}$, and $j_{\infty}(k;\tau)=1$;  for $x\in\bbR$ choose
$\tau_{x}=-1/(\tau-x)$,  $q_{x}=e^{2\pi\I\tau_{x}/h_x}$  and 
$j_{x}(k;\tau)=\tau_{x}^{k}$; while for $z\in\uhp$ choose $\tau_{z}=(\tau-z)/(\tau-\overline{z})$, 
$q_{z}=\tau_{z}^{n_z}$ and $j_{z}(k;\tau)=(1-\tau_{z})^{k}$. This factor $j_{z}$
is, up to a constant, the standard weight-$k$ automorphy factor associated to the transformation
$\tau\mapsto\tau_z$. 

The point is that any meromorphic function $f(\tau)$  invariant under the slash operator
\begin{equation}(f|_k\gamma_{z;h})(\tau):=(c\tau+d)^{-k}f\left(\frac{a\tau+b}{c\tau+d}\right)\,,\end{equation}
for some $z\in\bbH_\Gamma$, where we write $\gamma_{z;h}=\left({a\atop c}{b\atop d}
\right)$,   will have a local expansion
\begin{equation}f(\tau)=j_{z}(k;\tau)\,q_z^{k/n_z}\sum_{n\in\bbZ}f\left[n+\frac{k}{n_z}\right]_{z}q_{z}^{n}\,.
\label{eq:laurexp1}\end{equation}
The \textit{order} ord$_z(f)$ of an automorphic form $f$ at a point $z\in\bbH_\Gamma$ is defined to be the smallest 
$r\in\bbQ$ such that $f[r]_z\ne 0$. 

A \textit{quasi-automorphic form $f$} of weight $k\in 2\bbZ$ and depth $\le p$ for $\Gamma$ 
can be defined \cite{BGHZ} as a 
 function meromorphic on $\bbH_\Gamma$
(meromorphicity at the  cusps  is defined shortly), satisfying the functional equation
\begin{equation}
(f|_k\gamma)(\tau)=\sum_{r=0}^pf_r(\tau)\left(\frac{c}{c\tau+d}\right)^r\qquad\forall\gamma=\left({a\atop c}{b\atop d}\right)\in\Gamma
\end{equation}
for some functions $f_r$ meromorphic in $\bbH_\Gamma$ and independent of $\left({a\atop c}{b\atop d}\right)$.  We say $f$ is \textit{meromorphic at the cusp} $x\in\{\I\infty\}\cup\bbR$ if  all but finitely many coefficients $f[n]_x$ vanish
for $n<0$, and \textit{holomorphic at} $x$ if $f[n]_x=0$ whenever the relevant power of
$q_x$, namely $n+k/h_x$, is negative. 
When $p=0$, $f$ is called an \textit{automorphic form}; when $p=k=0$, it is called an \textit{automorphic function}.
When $\Gamma$ is commensurable with $\Gamma(1)$, it is typical to replace `automorphic' with
`modular'.

This definition can be extended to any weight
$k\in\bbC$ using the notion of automorphy factor, but we don't need it (though see the
end of Section 2.4). 
It is elementary to verify that the orders ord$_z(f)$ of an automorphic form $f$ are constant
on $\Gamma$-orbits $\Gamma z$. 

Suppose $f$ is an automorphic function, not constant. Then $f'= \frac{\mathrm{d}}{\mathrm{d}\tau}
f$ will be an automorphic form of weight 2, and $e_{2,f}=\frac{1}{f'}\frac{\mathrm{d}^2}{\mathrm{d}\tau^2}f$
will be quasi-automorphic of weight 2 and depth $1$. In this case, the \textit{Serre derivative} $D_k=\frac{\mathrm{d}}{\mathrm{d}\tau}-{k}{\beta}e_{2,f}(\tau)$,
 for some constant $\beta\in\bbC$ independent of $f$ and $k$ (computed for triangle groups in Theorem \ref{lemma1}(ii)
below),
takes automorphic forms of weight $k$ to those of weight $k+2$.
 
The automorphic functions form a field; when the genus of $\Gamma\backslash\bbH_\Gamma$ is zero,
this field can be expressed as the rational functions $\bbC(f)$ in some generator $f$. By a \textit{Hauptmodul} we
 mean any such generator. These Hauptmoduls $f$
are mapped to each other by the M\"obius transformations PSL$(2,\bbC)$, and therefore are
determined by 3 complex parameters. 

For example, for $\Gamma(1)=\mathrm{PSL}(2,\bbZ)$, recall
 the classical Eisenstein series $E_k$ given by 
\begin{equation}\label{eisen}
E_k(\tau)=1+\sum_{n=1}^\infty \frac{n^{k-1}q^n}{1-q^n}\,,\end{equation}
$k\in\bbZ_{>0}$, where $q=q_{\I\infty}=
\exp(2\pi\I\tau)$. The holomorphic modular forms and quasi-modular forms yield the polynomial 
rings $\bbC[E_4,E_6]$ and $\bbC[E_2,E_4,E_6]$.
The classical Hauptmodul is 
\begin{equation}\label{jsl2Z}j(\tau)=\dfrac{1728 E_4(\tau)^3}{E_4(\tau)^3-E_6(\tau)^2}=q^{-1}+744+196884q+\cdots\,.\end{equation} 
Throughout this paper, by $E_k
(\tau)$ and $j(\tau)$ we mean
these modular forms for $\Gamma(1)$.

\subsection{Triangle groups}
In this paper we focus on the \textit{triangle groups}. These by definition are
those genus-0 Fuchsian groups $\Gamma$ of the first kind with $n_{el}+n_{cp}=3$ (the minimal value possible). This means there are exactly 3 $\Gamma$-orbits of cusps and elliptic
fixed-points, in some combination. Let $2\le m_1\le m_2\le m_3\le \infty$ be the orders of the 
stabilizers of those 3 orbits. No Fuchsian group of the first kind can have 
 types $(2,2,m)$ $\forall m\le \infty$, $(2,3,n)$ for $n\le 6$, (2,4,4) and
(3,3,3); the remainder are called the hyperbolic types.  
We are primarily interested in the case where  $m_3=\infty$ --- for $m_3<\infty$ see 
Appendix B.
 As an abstract group, a triangle group has presentation $\langle g_1,g_2,g_3\,|\,g_i^{m_i}=1=g_1g_2g_3\rangle$; when $m_3=\infty$ this is isomorphic to the free
product $\bbZ_{m_1}*\bbZ_{m_2}$, where we write $\bbZ_k$ for the cyclic group with $k$ elements.

Given one such triangle group, we can find another by conjugating by any $g\in\mathrm{PSL}(2,\bbR)$.
The triangle group of a given type $\mft=(m_1,m_2,\infty)$ is unique up to this conjugation \cite{Pet}, and so is determined by 3 real
parameters. As the automorphic functions of $\Gamma$ and $g\Gamma g^{-1}$ are
related by $f(\tau)\leftrightarrow f(g^{-1}\tau)$, it is not so significant which realization is
chosen.  Of course,
this conjugation will in general affect the integrality of Fourier coefficients, so in that sense
some choices are better than others.

Write  $v_i=1/m_i$ for the angular parameters. 
A fundamental domain for a triangle group will be the double of a hyperbolic triangle in $\uhp_\mft$; we
fix the triangle group by fixing  the location of the corners of the triangle, which we take to be
\begin{equation}\label{corners}
\zeta_1=-e^{-\pi\I v_1}\,,\ \zeta_2=e^{\pi\I v_2}\,,\ 
\zeta_3=\I\infty
\,.\end{equation}
The Fuchsian group $\Gamma_{\mft}$ for this choice has  generators 
\begin{equation}\label{trigen}\gamma_1=\left(\begin{matrix}2\cos(\pi v_1)&1\cr -1&0\end{matrix}\right)\,,
\gamma_2=\left(\begin{matrix}0&1\cr -1&2\cos(\pi v_2)\end{matrix}\right)\,,
\gamma_3=\left(\begin{matrix}1&2\cos(\pi v_1)+2\cos(\pi v_2)\cr 0&1\end{matrix}\right)
\end{equation}
stabilizing the 3 corners $\zeta_1,\zeta_2,\zeta_3$, where 
\begin{equation}\gamma_1\gamma_2\gamma_3=\gamma_1^{m_1}=\gamma_2^{m_2}=-I_{2\times 2}\,.\end{equation}
Thus the cusp $\I\infty$ has cusp-width $h_{3}:=2\cos(\pi v_1)+2\cos(\pi v_2)$; when $m_2=\infty$,
 $\zeta_2=1$ is also a cusp, with cusp-width $h_2=1$.
Of course the groups $\Gamma_{(m_{\pi 1},m_{\pi 2},m_{\pi 3})}$  are conjugate for any permutation $\pi\in
\mathrm{Sym}(3)$.

The prototypical example is the modular group $\Gamma_{(2,3,\infty)}=\Gamma(1)$.
More generally, the \textit{Hecke groups} $\Gamma_{(2,m,\infty)}$, $m>2$, have attracted
a fair amount of attention.

\subsection{A Hauptmodul for triangle groups}


Given a  type ${\mft}={(m_1,m_2,\infty)}$, fix the triangle group  $\Gamma_{\mft}$ as in \eqref{trigen}. 
A  Hauptmodul $J_\mft(\tau)$ for $\Gamma_\mft$ is determined by 3 independent {\it complex} parameters,
which we fix  by  demanding
\begin{equation} \label{Jnorm}
J_{\mft}(\zeta_1)=1\,,\ J_{\mft}(\zeta_2)=0\,,\ J_{\mft}(\I\infty)=\infty\,.\end{equation}
(We make this choice because $1728J_{(2,3,\infty)}$ then equals the classical
choice \eqref{jsl2Z} for $\Gamma(1)$.)
We call the unique Hauptmodul satisfying \eqref{Jnorm} the \textit{normalized Hauptmodul} for $\Gamma_{\mft}$. To find it, given any other 
Hauptmodul $J$, first note that $J(\zeta_i)$ must be distinct points in
$\bbC\bbP^1$ (since $J$ is a Hauptmodul) so there will be a unique M\"obius transformation mapping those 3 points
to $1,0,\infty$ respectively, and $J_\mft$ is the composition of that transformation with $J$. Note that $J_{(m_1,\infty,m_2)}=J^{-1}_{(m_1,m_2,\infty)}$,
 $J_{(m_2,\infty,m_1)}=(1-J_{(m_1,m_2,\infty)})^{-1}$, etc.
In the following theorem we explicitly compute $J_{\mft}$, and in the
following section do this in a different way.

\begin{theo} \label{first} Fix any hyperbolic type $\mft=(m_1,m_2,\infty)$, $m_1\le m_2\le\infty$. 
Let ${q}_i$ be the local coordinates about the points $\zeta_i\in\bbH_\mft$ in \eqref{corners}, and
write $\widetilde{q}_i=\alpha_iq_i$ for $\alpha_i$ defined by: if $m_i=\infty$, 
\begin{equation}
\label{a-3}\alpha_i={b'd'}\prod_{k=1}^{b'-1}(2-2\cos(2\pi \frac{k}{b'}))^{-\frac{1}{2}\cos(2\pi \frac{ka'}{b'})}\prod_{l=1}^{d'-1}(2-2\cos(2\pi \frac{l}{d'}))^{-\frac{1}{2}\cos(2\pi \frac{lc'}{d'})}\,,\end{equation}
where we define integers $a',b',c',d'$ by $a'/b'=(1+v_1-v_2)/2$ and $c'/d'=(1+v_1+v_2)/2$; 
if $m_i<\infty$, 
\begin{equation}\label{a-1}\alpha_i=\frac{\cos(\pi (v_1+v_2)/2)}{\cos(\pi (v_1-v_2)/2)}
\frac{\Gamma(1+v_i)\Gamma((1-v_i+v_{3-i})/2)^2}{\Gamma(1-v_i)\Gamma((1+v_1+v_2)/2)^2}\,. 
\end{equation}
The normalized $J_{\mft}$ in \eqref{Jnorm} has local expansions
\begin{equation}J_{\mft}(\tau)=1+\widetilde{q}_1+\sum_{k=2}^\infty a_k\widetilde{q}_1^k=\widetilde{q}_2+\sum_{k=2}^\infty b_k\widetilde{q}_2^k=\widetilde{q}_3^{-1}+
\sum_{k=0}^\infty c_k\widetilde{q}_3^k\,,\label{jmnqexp}\end{equation}
These (normalized) coefficients $a_k,b_k,c_k$ are uniquely determined by
\begin{equation}\label{schw}-2{\dddot{J}_{\mft}\,\dot{J}_\mft}+{3}{\ddot{J}_{\mft}^2}
-n_z^{-2}\dot{J}_\mft^2=\dot{J}_{\mft}^4\left({1-v_2^2\over J_{\mft}^2}+{1-v_1^2\over (J_{\mft}-1)^2}+{v_1^2+v_2^2-1\over J_{\mft}(J_{\mft}-1)}\right)\,\end{equation}
for the choice $z=\zeta_1,\zeta_2,\zeta_3$ respectively, 
where each dot denotes $\widetilde{q}_j\frac{{\rm d}}{{\rm d}\widetilde{q}_j}$, and where $n_z$
is the order of the stabilizer at $z$. 
The coefficients $a_k,b_k,c_k$ are universal (i.e. type-independent) polynomials in $\bbQ[v_1,v_2]$, and are also unchanged if we replace $\Gamma_\mft$ by any conjugate.
 \end{theo}
 
 The key to this calculation, which we describe in Section 7.1, is the expression (using ratios of 
 hypergeometric functions) of the uniformizing Schwarz map from the upper hemisphere in $\bbC\bbP^1$
 to a hyperbolic triangle in the Poincar\'e disc. Analytically continuing the (multivalued) hypergeometric
 functions amounts to reflecting in the sides of that triangle, resulting in a multivalued
 map from the thrice-punctured sphere to the disc. The (single-valued) functional inverse of this
 Schwarz map is a Hauptmodul; its automorphy traces back to the monodromy
 of the hypergeometric equation. The most convenient way to obtain (most of) the local
 expansion of that Hauptmodul is through the Schwarzian equation \eqref{schw}.

For instance we have
\begin{align}
c_0&\,=(1-\gamma_-)/2\,,\ 
c_1=(5-2\gamma_+-3\gamma_-^2)/64\,,\ 
c_2=(-\gamma_-^3-\gamma_+\gamma_-+2\gamma_-)/54\,,\cr
c_3&\,=(-31+76\gamma_+-28\gamma_+^2+690\gamma_-^2-404\gamma_+\gamma_-^2-303\gamma_-^4)/32768\,,\cr 
c_4&\,=(-274\gamma_-+765\gamma_+\gamma_--314\gamma_+^2\gamma_-+2807\gamma_-^3-1865\gamma_+\gamma_-^3-1119\gamma_-^5)/216000\,,\cr
c_5&\,=(19683- 121770\gamma_++199044\gamma_+^2  - 1909439\gamma_-^2+5990732\gamma_+\gamma_-^2-68472\gamma_+^3 \cr &\ \ +12854105\gamma_-^4-2699804\gamma_+^2\gamma_-^2-9509386\gamma_+\gamma_-^4-4754693\gamma_-^6)/1528823808\,,\cr
c_6&\,=(341510 \gamma_-- 2360379 \gamma_+\gamma_-- 13805911 \gamma_-^3 + 4269300 \gamma_+^2\gamma_-- 1587244 \gamma_+^3\gamma_-+ 48264782 \gamma_+\gamma_-^3\cr&\ \ + 70933968 \gamma_-^5- 23644656 \gamma_+^2\gamma_-^3-57687959 \gamma_+ \gamma_-^5-24723411\gamma_-^7)/12644352000\,.
\nonumber\end{align}
where $\gamma_{\pm}=v_1^2\pm v_2^2$.  To our knowledge, these formulas in this generality have not appeared in the literature, although
\cite{wol83} computed \eqref{a-3}-\eqref{a-1}. Replacing $\Gamma_\mft$ with any conjugate 
(a 3 \textit{real} number ambiguity coming from PSL$(2,\bbR)$) affects $J_{\mft}$ by changing 
 the value of $\alpha_3$, the value of cusp-width $h_3$, and the choice of $\I\infty$ as a 
cusp. The only subtlety here is which $\alpha_3$ corresponds to our choice \eqref{trigen} of
$\Gamma_\mft$. We find that once one has chosen $\I\infty$ to be a cusp (it could have been 
anywhere in $\bbR\cup\{\I\infty\}$) and has fixed
the cusp-width $h_3$ (it could have been any positive real number), then the modulus $|\alpha_3|$ is
fixed for any conjugate; our choice \eqref{trigen} of generators then corresponds to $\alpha_3$ being positive.

\subsection{Automorphic forms for triangle groups} 

Knowing a Hauptmodul $J$ for any genus-0 Fuchsian group --- e.g. any triangle group ---
determines by definition all automorphic functions. It is less well known that from a Hauptmodul,
all holomorphic (quasi-)automorphic forms can be quickly read off.  
 We restrict here to triangle groups, although the argument works for
any genus-0 group.

The following theorem constructs an automorphic form whose divisor is
supported at the cusps, the analogue here of the discriminant
form $\Delta=\eta^{24}$ for $\Gamma(1)$. It constructs from this a `rational' basis for
the space of automorphic forms (rational in a sense described after the theorem), and  gives
the analogue here of $E_2$, and hence all quasi-automorphic forms. In Section \ref{arith11apr}, we
compare this basis with more classical ones, for the 9 triangle groups related to $\Gamma(1)$.

\begin{theo}\label{lemma1} \begin{enumerate}
\item[(i)] {For each $k\in \bbZ$, write $d_{2k}=k-\lceil {k/m_1}\rceil-\lceil{k/ m_2}\rceil$
and let 
\begin{equation}f_{2k}=({-1})^k\dot{J}_\mft^kJ_\mft^{\lceil \frac{k}{m_2}\rceil-{k}}(J_\mft-1)^{\lceil \frac{k}{m_1}\rceil-{k}}=\widetilde{q}_3^{d_{2k}}+O(\widetilde{q}_3^{d_{2k}+1})\,,
\label{modforms}\end{equation}
where the dot denotes $\widetilde{q}_3\mathrm{d}/\mathrm{d}\widetilde{q}_3$.
Then a basis for 
the $\C$-vector space $\frak{m}_{2k}(\Gamma_\mft)$ of holomorphic automorphic forms of 
weight $2k$ for $\Gamma_\mft$ is  $f_{2k}(\tau)\,J_{\mft}(\tau)^l$ for each $0\le l\le d_{2k}$. In particular,} \begin{equation}
\mathrm{dim}(\frak{m}_{2k}(\Gamma_\mft))=\left\{\begin{matrix}d_{2k}+1&\mathrm{if}\ k\ge 0\\
0&\mathrm{if}\ k<0\end{matrix}\right.\,.\end{equation} 

The algebra $\frak{m}(\Gamma_\mft)$ of holomorphic automorphic forms has the following minimal set of  generators:

when $\mft=(\infty,\infty,\infty)$, $\{f_2,J_\mft f_2\}$;

when $\mft=(m,\infty,\infty)$ for $m<\infty$, $\{f_2,\ldots,f_{2m}\}$;

when $\mft=(m_1,m_2,\infty)$ for $m_1\le m_2<\infty$, $\{f_{2l}\}|_{2\le l\le m_2}\cup\{J_\mft^{d_{2l}}f_{2l}\}|_{3\le l\le m_1}$\,.

\item[(ii)] {Define $L$ to be the least common multiple $\mathrm{lcm}(m_1,m_2)$ where we 
write lcm$(m_1,\infty)=m_1$ and lcm$(\infty,\infty)=1$. Then 
$\Delta_{\mft}(\tau):=f_{2L}(\tau)$ 
is a holomorphic automorphic form of weight $2L$, nonzero 
everywhere in $\bbH_\mft$ except in the $\Gamma_\mft$-orbit $[\I\infty]$, where $\Delta_\mft$ has a
zero of order $n_\Delta=L\,(1
-m_1^{-1}-m_2^{-1})$. Define $E_{2;\mft}=\frac{1}{2\pi\I}\Delta_\mft^{-1}\mathrm{d}{\Delta_\mft}/\mathrm{d}\tau$.
Then $E_{2;\mft}$ is holomorphic in $\uhp_\mft$, $E_{2;\mft}$ vanishes at any cusp
$\zeta_{m_j}\not\in[\I\infty]$, and $E_{2;\mft}(\I\infty)=n_\Delta$. Moreover, $E_{2;\mft}$ is quasi-automorphic of weight 2 and depth 1  for $\Gamma_\mft$: i.e. for
all $\left({a\atop c}{b\atop d}\right)\in\Gamma_\mft$, 
\begin{equation}
E_{2;\mft}\left({a\tau+b\over c\tau+d}\right)={n_\Delta c\over 2\pi\I}(c\tau+d)E_{2;\mft}(\tau)+
(c\tau+d)^2E_{2;\mft}(\tau)\,.
\end{equation}
The derivation 
$$D_k=\frac{1}{2\pi\I}\frac{\mathrm{d}}{\mathrm{d}\tau}-\frac{k}{L}E_{2;\mft}$$
sends weight $k$ automorphic forms to weight $k+2$ ones.
The space of all holomorphic quasi-automorphic forms of $\Gamma_\mft$ is $\frak{m}(\Gamma_\mft)[E_{2;\mft}]$.} 
\end{enumerate}\end{theo}    
   
   The $f_{2k}$ defined above is the unique normalized holomorphic weight-$2k$ automorphic
   form with maximal order at the cusp $\I\infty$.  
The weights of generators for $\frak{m}(\Gamma)$ for any Fuchsian group of the first kind, 
are given in \cite{wag} and references therein; what we provide in Theorems 1 and 2 are
explicit  formulas and expansions for those generators, in the special case of triangle groups. 
Provided we expand  in $\widetilde{q}_i=\alpha_iq_i$ instead of
$q_i$, $J_{\mft}$ has rational coefficients; in this same sense, our bases for each $\mathfrak{m}_{2k}$ also has rational coefficients. Incidentally, according to
Wolfart \cite{wol83}, $\alpha_3$ is 
transcendental except for the types listed in Table 1 below. 

Although every triangle group shares many properties with $\Gamma(1)$, one difference is that 
$\frak{m}(\Gamma_\mft)$
will rarely be a polynomial algebra: in fact, $\frak{m}(\Gamma_\mft)$ is polynomial iff
$\mft=(2,3,\infty),(2,\infty,\infty)$, or $(\infty,\infty,\infty)$.
On the other hand, \cite{Mil,wol81} consider the ring of holomorphic automorphic forms of $\Gamma_\mft$ for 
a root-of-unity-valued multiplier (which allows certain weights $k\not\in 2\bbZ$), and find that that
larger ring always generated by 3 forms $f_1,f_2,f_3$ satisfying an identity of the form 
$f_1^{e_1}+f_2^{e_2}+
f_3^{e_3}=0$.

Incidentally,  $\Delta_\mft$ can identify all automorphic forms with
 multiplier of arbitrary \textit{complex} weight $k\in\bbC$. In particular, for any $w\in\bbC$
 define $\Delta^{(w)}_\mft$ to be any nontrivial  solution to
\begin{equation}\label{eq:ODEDelta}
\frac{1}{2\pi\I}\frac{\mathrm{d}}{\mathrm{d}\tau}{f}=wE_{2;\mft}f\,.\end{equation} 
First note from the theory of ordinary differential equations (see e.g. \cite{bib:Hi}), $\Delta_\mft^{(w)}$ exists and is
holomorphic throughout $\uhp$. Locally, it corresponds to some branch of the power $\Delta_\mft^w$;
that it transforms under $\Gamma_\mft$ like (and therefore is) a holomorphic automorphic form of weight $w\,\mathrm{lcm}
\{m_1.m_2\}$ follows directly from \eqref{eq:ODEDelta}. Then some $f$ is a (meromorphic)
automorphic \textit{form} for $\Gamma_\mft$ with arbitrary  weight $k\in\bbC$ automorphy factor, iff 
$f/\Delta_\mft^{(k/\mathrm{lcm}\{m_1,m_2\})}$ is an automorphic \textit{function} for
$\Gamma_\mft$ with the appropriate automorphy factor (namely some character of
$\Gamma_\mft$).

\section{Quasi-automorphic forms via Halphen's equation}
\label{24oct2012}
In this section we realize the (quasi-)automorphic forms of the triangle groups, using the Halphen
differential equation. This material should be completely new;
see \cite{ho07-1} for some of the detailed calculations which are omitted here. 
 For simplicity, we again require 
$m_3=\infty$ --- see Appendix B for some remarks on the generalization to finite $m_3$. 


Fix any hyperbolic type $\mft=(m_1,m_2,\infty)$. 
 Recall the angular parameters $v_i=1/m_i$. 
 Consider the Halphen differential equation \eqref{tatlis},
where $a,b,c$ are the parameters 
\begin{eqnarray*}
a &=&\frac{1}{2}(1+v_2-v_1-v_3)\,,\\
b &=& \frac{1}{2}(1+v_3-v_1-v_2)\,,\\
c &=& \frac{1}{2}(1+v_1-v_2-v_3)\,.
\end{eqnarray*}
In the original Halphen  equation, the right hand side of  (\ref{tatlis}) is divided by $a+b+c-2$.

Recall the 
normalized Hauptmodul $J_\mft$.
We are interested in the particular solution of \eqref{tatlis} given in Theorem 3(i) below.
Because $v_3=0$ (i.e. $a+c=1$), the Halphen vector field has the one-dimensional singular locus 
$t_1=t_3=0$; the solution of part (i) is a perturbation of this singular locus.
The relation of the Halphen equation with hypergeometric functions goes back to Halphen,
who is therefore ultimately  responsible for part (i), (iii).    
Part (ii) follows from recursions coming from \eqref{tatlis} (see Section 7.3 below), and is new.
The automorphy of the Halphen solutions arises from the $\SL2\C$ action in part (iii), and
can be also proved using generalizations of  period maps, see Section 10 of \cite{ho07-1}.

\begin{theo}
\label{15may2013-1}
\begin{enumerate}
\item[(i)]\label{9feb2012}
A solution to \eqref{tatlis}
is: 
\begin{eqnarray*}
 t_1(\tau)&=&  (a-1)z\,Q(z)\, F(1-a,b,1; z)\,F(2-a,b,2; z)\,,\notag\\
\label{t123}
t_2(\tau)&=&  Q(z)\,F(1-a,b,1; z)^2+t_1(\tau)\,,\\
t_3(\tau)&=&  Q(z)\,z\,F(1-a,b,1; z)^2+t_1(\tau)
\,,\notag\end{eqnarray*}
where $F={}_2F_1$ is the hypergeometric function and
$$
Q(z)=\frac{\pi\I\,(1-b)}{2 \sin(\pi b)\,\sin(\pi a)}(1-z)^{b-a}\, ,  \   \  z=(1-J_\mft(\tau))^{-1}
$$
\item[(ii)] Write $\qh=\nu e^{2\pi\I\tau/h_3}$ where $ h_3 =2\cos(\pi v_1)+2\cos (\pi v_2)$ and
\begin{equation}
\label{nuoct20}
\nu =  \left \{ \begin{array}{ll}
          \frac{1}{2}v_1^2v_2^2 \alpha_3 & v_1\not=0\,,\ v_2\not=0\,,\\
              \frac{1}{2}v_1^2 \alpha_3 & v_2=0\,,\ v_1\not=0\,, \\ 
               8& v_1=0,\ v_2=0\,. 
          \end{array} \right. 
\end{equation}
Then the solution of (i) has the expansion
\begin{equation}\label{ti}
t_i=\frac{2\pi \I}{h_3} t_{i,0}+\kappa_i\sum_{j=1}^\infty \tilde{t}_{i,j}\qh^j,\ 
\end{equation}
where $[t_{1,0},t_{2,0},t_{3,0}]=[0,-1,0]\,,$ and
$$
[\kappa_1,\kappa_2,\kappa_3]=
\frac{2\pi\I}{h_3}\left[
\begin{array}{*{3}{c}}
-m_{1}^{2}m_2^{2}-m_{2}^{2}m_1+m_{2}m_1^{2}, & m_2m_1+m_2+m_1, & 
m_{1}^{2}m_2^{2}-m_{2}^{2}m_1+m_{2}m_1^{2}
\end{array}
\right]\,,
$$
\begin{equation}
\label{31may2012}
\tilde{t}_{i,j}\in\Q[m_1,m_2]\,. 
\end{equation}
\item[(iii)]
\label{1feb11}
 If $t_i(\tau),\ i=1,2,3,$  are the coordinates of 
any solution of the Halphen differential equation,  then so are 
$$
\frac{1}{(c'\tau+d')^2}\,t_i\left(\frac{a'\tau+b'}{c'\tau+d'}\right)-\frac{c'}{c'\tau+d'},\ \ \forall\mat{a'}{b'}{c'}{d'}\in\SL 2\C.
$$
\end{enumerate}
\end{theo}\medskip
For example, $\tilde{t}_{1,1}=\tilde{t}_{3,1}=1$, $\tilde{t}_{2,1}=m_1-m_2$,
$$\tilde{t}_{1,2}=\frac{1}{4}(2m_1m_2^2-m_1^2m_2^2-7m_1^2+7m_2^2)\ ,\ 
\tilde{t}_{3,2}=\frac{1}{4}(m_1^2m_2^2-7m_1^2+7m_2^2-2m_1^2m_2)\,,$$
$$\tilde{t}_{2,2}=\frac{1}{8}(-m_1^3m_2^3+6m_1^2m_2^2-11m_1^3+11m_1^2m_2-m_1^3m_2^2-3m_1^3m_2-11m_2^3-m_1^2m_2^3+11m_1m_2^2-3m_1m_2^3)\,,$$
$$\tilde{t}_{1,3}=\frac{1}{48}(3 m_1^4  m_2^4-14m_1^2  m_2^4  - 64 m_1^3  m_2^2    + 64 m_1 m_2^4
+50m_1^4m_2^2+139m_1^4+139m_2^4-278m_1^2m_2^2)\,,$$
$$\tilde{t}_{3,3}=\frac{1}{48}(3m_1^4m_2^4-14m_1^4m_2^2+64m_1^4m_2+139m_1^4-64m_1^2m_2^3+139m_2^4-278m_1^2m_2^2+50m_1^2m_2^4)\,.$$

Recall the triangle group $\Gamma_\mft$ of type $\mft=(m_1,m_2,\infty)$
generated by the matrices \eqref{trigen}. We focus in this section on 
$\hat{q}$-expansions around the cusp $\I\infty$. 
The renormalization by $\nu$ of $\alpha_3$ is natural from the point of view of the
recursion coming from \eqref{tatlis}. 
For each $ k\ge 2$, we set
\begin{eqnarray*} 
E_{2k,\mft}^{(1)} &:=& \left(\frac{h_3}{2\pi\I}\right)^k(t_1-t_2)\,( t_3-t_2)^{k-1}\in 1+\hat{q}\,\bbQ[[\hat{q}]]\,,  \\
E_{2k,\mft}^{(2)} &:= & \left(\frac{h_3}{2\pi\I}\right)^k(t_1-t_2)^{k-1}\,( t_3-t_2)\in 1+\hat{q}\,\bbQ[[\hat{q}]]\,,\\
E_{4,\mft}&:=& E_{4,\mft}^{(1)}=E_{4,\mft}^{(2)}\,, \\
E_{6,\mft} &:=& E_{6,\mft}^{(2)}\,. 
\end{eqnarray*}
Define $E_{2,\mft}$ using Theorem 4(iii).
 The notation and normalization is chosen so that when $\mft=(2,3,\infty)$,
$E_{k,\mft}$ for $k=4,6$ coincide  with the 
classical series for $\Gamma(1)$. From now on we regard all $t_i$'s as functions of $\tau$.
The convention throughout this paper is that the value of a polynomial $P(x)$ for $x=\infty$ is the coefficient of the  monomial $x^n$ of highest degree in $P(x)$.

\begin{theo}
\label{29feb2012}
Assume as usual that $2\le m_1\le m_2\le\infty$ and $\mft=(m_1,m_2,\infty)$ is hyperbolic. Then
\begin{enumerate}
\item[(i)] The $t_i(\tau)$ are quasi-automorphic. More precisely,
 they are meromorphic functions of $\tau\in\uhp_\mft$, and
satisfy the following functional equation:
\begin{equation}
\label{aruquia}
(c'\tau+d')^{-2}t_i(\gamma(\tau))-c'(c'\tau+d')^{-1}=t_i(\tau)\qquad\forall\gamma=\left({a'\atop c'}
{b'\atop d'}\right)\in\Gamma_\mft\,.
\end{equation} 
 \item[(ii)]
The field generated by all meromorphic automorphic forms 
for $\Gamma_\mft$ consists of all rational functions in $t_1-t_2$ and $t_3-t_2$. 
\item[(iii)] The relation with Theorems 1 and 2 is:   $t_1-t_2=\frac{2\pi\I}{h_3}\frac{J'_\mft}{J_\mft}$ and $t_3-t_2=\frac{2\pi\I}{h_3}\frac{J'_\mft}{(J_\mft-1)}$,
\begin{eqnarray*}
\frac{1}{n_\Delta}{E}_{2;\mft} &=& 
\frac{b-a}{b}t_1-t_2+\frac{a+b-1}{b}t_3\,,  \\
f_4&=&E_{4,\mft}\,,\ \ \ f_6=\left\{\begin{matrix} E_{6,\mft}&\mathrm{if}\ m_1=2\\ E_{6,\mft}/(J_\mft-1)&\mathrm{otherwise}\end{matrix}\right.\,,\\
J_{\mft} &=& \frac{t_3-t_2}{t_3-t_1}=\frac{E_{4,\mft}^3}{E_{4,\mft}^3-  E_{6,\mft}^2}\,.
\end{eqnarray*}
Moreover, the function $j_\mft=2m_2^2m_1 ^2J_\mft+(-m_2^2m_1^2+m_2^2-m_1^2)$ is the unique Hauptmodul for $\Gamma_\mft$ normalized so that $j_\mft(\tau)=\frac{1}{\qh}+O(\qh^1)$. 
\item[(iv)] When $m_2\ne \infty$, the algebra $\mathfrak{m}(\Gamma_\mft)$ of holomorphic automorphic forms is generated by
$$
E_{2k,\mft}^{(2)}\,,\ 2\leq k \leq m_2\,,\ \ E_{2k,\mft}^{(1)}\,,\ 3\leq k \leq m_1\,.
$$
When $m_1<\infty=m_2$, $\mathfrak{m}(\Gamma_\mft)$ is generated by
$$
E_{2k,\mft}^{(1)}\,,\ 1\leq k \leq m_1\,.
$$
\end{enumerate}
\end{theo}
The case $m_1=m_2=m_3=\infty$ corresponds to the classical Darboux-Halphen differential equation, see \S\ref{takeuchi}.

It should be emphasized that, although ultimately the approaches in Sections 2 and 3 both
reduce to hypergeometric calculations, the approaches are independent in the sense
that their outputs (a Hauptmodul in \S2 compared with three quasi-automorphic forms in \S3) are 
different. Both approaches are complete in the sense that all (quasi-)automorphic forms
for the given triangle group $\Gamma_\mft$ can be obtained from their outputs by standard operations.






\section{The modular triangle groups}
\label{arith11apr}

By a \textit{modular triangle group} $\Gamma$ we mean a triangle group commensurable with $\Gamma(1)$
(i.e. $\Gamma\cap\Gamma(1)$ has finite index in both $\Gamma$ and $\Gamma(1)$). There are precisely 9
$\Gamma_{\mft}$ conjugate to a modular triangle group \cite{tak77}. Such Fuchsian groups are called
\textit{arithmetic} (the definition of arithmetic Fuchsian groups can be extended to the case
where there are no cusps, and \cite{tak77} also identifies these). 
 In this section we show how our expressions for modular forms recover
the classical ones in these 9 cases.

In Table 1 we list these 9 types, together
with one of the modular triangle groups which realizes it. We include the
basic data for that conjugate $g\Gamma_{\mft}g^{-1}$.
In the table and elsewhere, we write $\omega=e^{2\pi\I/6}$, $S=\left({0\atop 1}{-1\atop 0}\right)$,
$T=\left({1\atop 0}{1\atop 1}\right)$, and $U=\left({1\atop -1}{0\atop 1}\right)$. The matrix
$W_N={1\over\sqrt{N}}\left({0\atop -N}{1\atop 0}\right)$ is called a \textit{Fricke involution}.
As usual, $\Gamma(N)$ consists of all $A\in\Gamma(1)$ with $A\equiv \pm I$ (mod $N$),
$\Gamma_0(N)$ consists of all $A\in \Gamma(1)$ with entry $A_{2,1}$ divisible by $N$, and
$\Gamma_0^+(N):=\langle \Gamma_0(N),W_N\rangle$. Given any triangle group $\Gamma$ of type
$(2,n,\infty$), by $\Gamma^*$ we mean the subgroup generated by the squares $\gamma^2$ of all elements $\gamma\in
\Gamma$, together with any element in $\Gamma$ of order $n$; then $\Gamma^*$ has index 2 in
$\Gamma$, and is a triangle group of type $(n,n,\infty)$. Table 1 is largely taken from
\cite{BT}.

\medskip\centerline{{\it Table 1.} The triangle groups commensurable with  $\Gamma(1)$}
$$\tiny\vbox{\tabskip=0pt\offinterlineskip
  \def\tablerule{\noalign{\hrule}}
  \halign to 5.7in{\strut#
  \tabskip=0em plus1em &    
    \hfil#&\vrule#&\hfil#&\vrule#&    
\hfil#&\vrule#&\hfil#&\vrule#&\hfil#&\vrule#&\hfil#&\vrule#&\hfil#&\vrule#&\hfil#&\vrule#&\hfil#
&\vrule#&\hfil#&\vrule#
\tabskip=0pt\cr
&$(m_1,m_2,m_3)\,\,$&&\hfil$\,g\Gamma_{\mft}g^{-1}$\hfil&&\hfil$g$\hfil&&\hfil$\zeta_1$\hfil&&\hfil$\gamma_1$\hfil&&\hfil$\zeta_2$\hfil&&\hfil 
$\gamma_2$\hfil&&\hfil$\,\zeta_3\,$\hfil&&\hfil$\,\gamma_3\,$ \hfil&&\hfil$\alpha_3$\hfil\cr
\tablerule&\hfil$(2,3,\infty)$\hfil&&$\Gamma(1)$\hfil&&\hfil 1\hfil&&\hfil$\I$\hfil&&\hfil$S$\hfil&&\hfil$\omega$\hfil&&\hfil
$\left({0\atop -1}{1\atop 1}\right)$\hfil&&\hfil$\infty$\hfil&&\hfil$T$\hfil&&\hfil$1728$\hfil\cr
\tablerule&\hfil$(2,4,\infty)$\hfil&&\hfil$\Gamma^+_0(2)$\hfil&&$\left({2\atop 0}{0\atop 1}\right)$
&&\hfil$\I/\sqrt{2}$\hfil&&\hfil$W_2$\hfil&&\hfil${(-1+\I)/2}$\hfil&&\hfil
$\frac{1}{\sqrt{2}}\left({2\atop-2}{1\atop 0}\right)$\hfil&&\hfil$\infty$\hfil&&\hfil$T$\hfil&&\hfil$256$\hfil\cr
\tablerule&\hfil$(2,6,\infty)$\hfil&&\hfil$\Gamma^+_0(3)$\hfil&&$\left({3\atop 1}{0\atop 1}\right)$&&\hfil$\I/\sqrt{3}$\hfil&&\hfil$W_3$\hfil&&\hfil
$(-3+\I\sqrt{3})/6$\hfil&&\hfil$\frac{1}{\sqrt{3}}\left({3\atop -3}{1\atop 0}\right)$\hfil&&\hfil$\infty$\hfil&&\hfil$T$\hfil&&\hfil$108$\hfil\cr
\tablerule&\hfil$(2,\infty,\infty)$\hfil&&$\Gamma_0(2)$\hfil&&$\left({1\atop 0}{1\atop 2}\right)$&&\hfil$(1+\I)/2$\hfil&&\hfil$\left({1\atop 2}{-1\atop-1}\right)$\hfil&&\hfil$0$\hfil&&\hfil
$U^{2}$\hfil&&\hfil$\infty$\hfil&&\hfil$T$\hfil&&\hfil $64$\hfil\cr
\tablerule&\hfil$(3,3,\infty)$\hfil&&\hfil$\Gamma(1)^*$\hfil&&\hfil1\hfil&&\hfil$\omega^2$\hfil&&\hfil$\left({1\atop-1}{1\atop 0}\right)$\hfil&&\hfil$\omega$\hfil&&\hfil
$\left({0\atop -1}{1\atop 1}\right)$\hfil&&\hfil$\infty$\hfil&&\hfil$T^2$\hfil&&\hfil$48\sqrt{3}$\hfil\cr
\tablerule&\hfil$(3,\infty,\infty)$\hfil&&$\Gamma_0(3)$\hfil&&$\left({1\atop 0}{-1\atop 3}\right)$&&\hfil$(3+\I\sqrt{3})/6$\hfil&&\hfil$\left({1\atop 3}{-1\atop-2}\right)$\hfil&&\hfil$0$\hfil&&\hfil
$U^{3}$\hfil&&\hfil$\infty$\hfil&&\hfil$T$\hfil&&\hfil$27$\hfil\cr
\tablerule&\hfil$(4,4,\infty)$\hfil&&\hfil$\Gamma^+_0(2)^*$\hfil&&$\left({2\atop 1}{0\atop 1}\right)$
&&\hfil$(\I-1)/2$\hfil&&\hfil$\frac{1}{\sqrt{2}}\left({2\atop -2}{1\atop 0}\right)$\hfil&&\hfil$(1+\I)/2$\hfil&&\hfil
$\frac{1}{\sqrt{2}}\left({0\atop -2}{1\atop 2}\right)$\hfil&&\hfil$\infty$\hfil&&\hfil$T^2$\hfil&&\hfil$32$\hfil\cr
\tablerule&\hfil$(6,6,\infty)$\hfil&&\hfil$\Gamma^+_0(3)^*$\hfil&&$\left({3\atop 1}{0\atop 1}\right)$
&&\hfil$(-3+\I\sqrt{3})/6$\hfil&&\hfil
$\frac{1}{\sqrt{3}}\left({3\atop -3}{1\atop 0}\right)$\hfil&&\hfil$(3+\I\sqrt{3})/6$\hfil&&\hfil$\frac{1}{\sqrt{3}}\left({0\atop -3}{1\atop 3}\right)$\hfil&&\hfil$\infty$\hfil&&\hfil$T^2$\hfil&&\hfil$12\sqrt{3}$\hfil\cr
\tablerule&\hfil$(\infty,\infty,\infty)$\hfil&&$\Gamma(2)$\hfil&&\hfil$\left({1\atop 0}{1\atop 2}\right)$\hfil&&\hfil$0$\hfil&&\hfil$U^{2}$\hfil&&\hfil$1$\hfil&&\hfil
$\left({-1\atop -2}{2\atop 3}\right)$\hfil&&\hfil$\infty$\hfil&&\hfil$T^2$\hfil&&\hfil$16$\hfil\cr
\noalign{\smallskip}}}$$

In this section we recover explicitly the classical result that:

\begin{prop}
\label{13march} The algebra of holomorphic 
 modular forms for each modular triangle group  has a basis in $\bbZ[[Q]]$, where $Q$
is some rescaling of $q$ or $q^{1/2}$.
\end{prop}


Indeed, by Lemma 3 of \cite{Ga}, $1728J_{(2,3,\infty)}$, $256J_{(2,4,\infty)}$, $108J_{(2,6,\infty)}$,
$16J_{(\infty,\infty,\infty)}$, $64J_{(2,\infty,\infty)}$, and $27J_{(3,\infty,\infty)}$ all have integer $q$-
or $q^{1/2}$-coefficients (whichever is appropriate), and leading term $\pm q^{-1}$ or
$\pm q^{-1/2}$.  $144J_{(3,3,\infty)}$,
$32J_{(4,4,\infty)}$, and $36J_{(6,6,\infty)}$ have $q^{1/2}$-coefficients in the Eisenstein
$\bbZ[\omega]$ or Gaussian $\bbZ[\I]$ integers, but if $Q$ is chosen to be $\I q^{1/2}/\sqrt{3}$,
$\I q^{1/2}$ or $\I q^{1/2}/\sqrt{3}$, respectively, then these functions lie in $Q^{-1}+\bbZ[[Q]]$. This 
information is  enough to verify that the basis given in Theorem 2 has integer coefficients.
The exact rescaling of $q$ or $q^{1/2}$ depends on the choice of realization of $\Gamma_\mft$.

\subsection{Type $\mft_m=(2,m,\infty)$  for $m=3,4,6$} 
For type $\mft=(2,3,\infty)$, the triangle group $\Gamma_\mft$ is  the full modular group $\Gamma(1)=\mathrm{PSL}(2,\bbZ)$. 
Its algebra of holomorphic quasi-modular forms is generated by the classical Eisenstein series
$E_2,E_4,E_6$ in \eqref{eisen}. We have $D_0=-E_4^2E_6/\Delta$.
Their relation with the quasi-modular forms coming from the Halphen system are
$$
E_{2;\mft}=E_2\,,\quad E_{4,\mft}^{}=E_4\,, \quad E_{6,\mft}^{}=E_6\,,\quad
J_{\mft}=j/1728\,.
$$

More generally, for any Hecke group $\Gamma_{(2,m,\infty)}$ (any $m\ge 3$),
Eisenstein series $E_{k,\mft_m}(\tau)$ can be analogously defined (see e.g. Section 4 of
\cite{leo08}).
The spaces of holomorphic automorphic forms of weights 4 and 6 are both one-dimensional,
spanned by what we call $f_4(\tau)=E_{4,\mft_m}(\tau)=1+\cdots$ and $f_6(\tau)=E_{6,\mft_m}(\tau)=
1+\cdots$ respectively. The normalized Hauptmodul is
\begin{equation}
J_\mft(\tau)=\frac{f_4(\tau)^3}{f_4(\tau)^3-f_6(\tau)^2}\,,\end{equation}
in perfect analogy with $\Gamma(1)$.
In the special cases $m=2p=4,6$ we are interested in here, we determine from Section 4.3.2 of \cite{leo08} that for any $k\ge 2$,
\begin{equation}\label{eisenhecke}E_{2k,\mft_{2p}}(\tau)=(E_{2k}(\tau)+p^{k}E_{2k}(p\tau))/(p^{k}+1)\end{equation}
and we find
\begin{align}\nonumber
J_{(2,4,\infty)}&=\frac{1}{256}q^{-1}+\frac{13}{32}+\frac{1093}{64}q+376q^2+\frac{620001}{128}
q^3+41792q^4+\cdots\,,\\ \nonumber
J_{(2,6,\infty)}&=\frac{1}{108}q^{-1}+\frac{1}{3}+\frac{371}{36}q+\frac{3643}{54}q^2+\frac{20713}{36}q^3-34396q^4+\cdots\,.
\end{align}

\subsection{Type $(\infty,\infty,\infty)$}
\label{takeuchi}

The most natural realization of $\mft=(\infty,\infty,\infty)$ is as $\Gamma(2)$, which has cusps at
$\I\infty,0,1$. The local parameter at the infinite
cusp is ${q}^{1/2}=e^{\pi \I\tau}$ (the square-root of the parameter for $\Gamma(1)$). 
Recall the Jacobi theta functions 
$$
\left \{ \begin{array}{l}
\theta_2(\tau):=\sum_{n=-\infty}^\infty q^{\frac{1}{2}(n+\frac{1}{2})^2}
\\
\theta_3(\tau):=\sum_{n=-\infty}^\infty q^{\frac{1}{2}n^2}
\\
\theta_4(\tau):=\sum_{n=-\infty}^\infty (-1)^nq^{\frac{1}{2}n^2}
\end{array} \right.\,. 
$$
It is well-known that $\theta_2^4,\theta_3^4,$ and $\theta_4^4=\theta_3^4-\theta_2^4$ are
modular forms for $\Gamma(2)$ of weight 2, and that they generate the ring of holomorphic
modular forms.  A Hauptmodul is 
$$J_{(\infty,\infty,\infty)}(\tau)=\frac{\theta_3(\tau)^4}{\theta_2(\tau)^4}=\frac{1}{16}q^{-1/2} + \frac{1}{2} +\frac{5}{4}q^{1/2}-  \frac{31}{8} q^{3/2}+\frac{27}{2}q^{5/2}+\cdots\,,$$
which maps $\I\infty$ to $\infty$, cusp 0 to 1 and cusp 1 to $0$.  
The normalized quasi-modular form is $e_2=E_2/6$.

In 1878 G. Darboux studied the system of differential equations
\begin{equation}
 \label{darboux}
\left \{ \begin{array}{l}
\dot u_1+\dot u_2=2u_1u_2\\
\dot u_2+\dot u_3=2u_2u_3\\
\dot u_1+\dot u_3=2u_1u_3
\end{array} \right.,
\end{equation}
in connection with  triply orthogonal 
surfaces in $\R^3$. Later Halphen in \cite{ha81}
found a solution of (\ref{darboux}) in terms of theta series:
$$
u_1=2(\ln \theta_4(\tau))'\,,
u_2=2(\ln \theta_2(\tau))'\,,
u_3=2(\ln \theta_3(\tau))'\,.
$$
The differential equation (\ref{darboux}) after the change of variables $t_i:=-2 u_i$ turns to be (\ref{tatlis}). The relations between the series $t_i$
in \S\ref{24oct2012} and theta series are given by
$$
\frac{-1}{4}t_i(8q^{\frac{1}{2}})=2 q\frac{\mathrm{d}}{\mathrm{d} q}\ln \theta_{j_i}\,,\ 
$$
where $ (j_1,j_2,j_3)=(3,2,4)$. 

\subsection{Types  $\mft_m=(m, \infty, \infty)$, $m=2,3$}
It is well-known that a Hauptmodul for $\Gamma_0(N)$ when $N-1$ divides 24 is $J_{(N)}(\tau)= 
(\eta(\tau)/\eta(N\tau))^{24/(N-1)}$, which for $N=2,3$ rescales to the normalized Hauptmoduln
\begin{align*}
&&J_{(2,\infty,\infty)}(\tau)=-\frac{1}{64}q^{-1}+\frac{3}{8}-\frac{69}{16}q+32q^2-\frac{5601}{32}q^3+768q^4-\frac{23003}{8}q^5+\cdots\,,\nonumber\\
&&J_{(3,\infty,\infty)}(\tau)=-\frac{1}{27}q^{-1}+\frac{4}{9}-2q+\frac{76}{27}q^2+9q^3-44q^4+\frac{1384}{27}q^5+
\cdots\,.\end{align*}
For any $N$ (and in particular $N=2,3$), 
$$q \frac{\mathrm{d}}{ \mathrm{d} q} \log\left( \frac{\eta(\tau ) }{\eta(N\tau )} \right)= E_2 (\tau )-N 
E_2 (N \tau ) $$
is a holomorphic weight-2 modular form for $\Gamma_0(N)$. 
For $\Gamma_0(2)$, the algebra of holomorphic modular forms is generated by
$E_2(\tau)-2E_2(2\tau)$ and $E_4(\tau)$, while that for  
$\Gamma_0(3)$ is generated by $E_2(\tau)-3E_2(3\tau)$, $E_4(\tau)$ and $E_6(\tau)$.

\subsection{Type $\mft'_m=(m,m,\infty)$ for $m=3,4,6$} 

Write $\mft_m=(2,m,\infty)$ as before.
Recall from the beginning of this section that a Fuchsian group of type  $\mft_m'$ (for any
$m\ge 3$) can be 
chosen to be the index 2 subgroup  $\Gamma^*_{\mft_m}$ of the Hecke group $\Gamma_{\mft_m}$.
The normalized Hauptmodul for  any $\mft_m'$ is 
$$J_{(m,m,\infty)}(\tau)=\frac{1}{2}\left(\frac{E_{6,\mft_m}(\tau)}{\sqrt{E_{6,\mft_m}(\tau)^2-E_{4,\mft_m}(\tau)^3}}+1\right)\,,$$
where $E_{k,\mft_m}=f_k$ here are the (normalized) Eisenstein series discussed in \S4.1.
The holomorphic  modular forms are generated by $\sqrt{E^2_{6,\mft_m}-E_{4,\mft_m}^3}$ together
with those for $\mft_m$ (since $\Gamma_{\mft_m'}$ is a subgroup of $\Gamma_{\mft_m}$). From this point of view the only thing special about
$m=3,4,6$ is that we can easily express $E_{4,\mft_m},E_{6,\mft_m}$ in terms of classical modular forms,
as was done in \eqref{eisenhecke} above. We find 
\begin{align}\nonumber
J_{(3,3,\infty)}(\tau)=&-\frac{\I\,\sqrt{3}}{144}q^{-1/2}+\frac{1}{2}+\frac{41\,\I\sqrt{3}}{12}q^{1/2}+
\frac{1255\,\I\sqrt{3}}{8}q^{3/2}+\frac{45925\,\I\sqrt{3}}{18}q^{5/2}+\cdots\,,\\ \nonumber
J_{(4,4,\infty)}(\tau)=&-\frac{\I}{32}q^{-1/2}+\frac{1}{2}+\frac{19\,\I}{8}q^{1/2}+\frac{351\,\I}{16}q^{3/2}+
\frac{653\,\I}{4}q^{5/2}+\frac{23425\,\I}{32}q^{7/2}+\cdots\,,\\ \nonumber
J_{(6,6,\infty)}(\tau)=&\frac{-\I\,\sqrt{3}}{36}q^{-1/2}+\frac{1}{2}+\frac{11\,\I\,\sqrt{3}}{12}q^{1/2}+\frac{17\,\I\,\sqrt{3}}{4} q^{3/2}+\frac{713\,\I\,\sqrt{3}}{36}q^{5/2}+\cdots\,.
\end{align}

\section{Observations and conjectures concerning coefficients}
The raison d'\^etre of modular forms is their $q$-expansions, i.e. the local (Fourier) expansions
about the cusp $\I\infty$. Expansions about other cusps have the same familiar feel (although
are usually ignored). The avoidance of considerations of 
(Taylor) expansions at points in $\bbH$, in particular at the elliptic fixed-points,
is almost complete.

It is hard to justify this focus on the
expansion at $\I\infty$, other than that it is exceedingly rich. However, a triangle group say has \textit{three}
 special $\Gamma_\mft$-orbits, perhaps the other two may also prove interesting. 
For example, in the vector-valued automorphic forms of Section 6.3 below, it seems artificial
to expand only about the large complex structure point (which corresponds to a cusp)
but to refuse to expand about say the Landau-Ginzburg point (which corresponds to an elliptic
fixed-point). For another example, consider the characters $\chi_M(\tau)=\sum_{r}a(M)_rq^{r}$
 of irreducible modules $M$ of rational vertex operator 
algebras. These $\chi_M$s  are modular functions for some $\Gamma(N)$. A surprise
happens at their expansions
$\chi_M(\tau)=\sum_ra(M)_{x;r}q_x^r$ about certain cusps $x\in\bbQ$ (which $x$ to choose
depends only on $N$): there are signs $\epsilon_x(M)$ and \textit{another} irreducible
module $M^x$ such that the coefficients at $x$ of $\chi_M$ equal those at $\I\infty$  of
$\epsilon_x(M)\,\chi_{M^x}$, that is, $a(M)_{x;r}=\epsilon_x(M)\,a(M^x)_r$. In other words,
expanding one character about a different cusp can recover a different character at
the usual cusp $\I\infty$. (This property of vertex operator algebra characters is implicit in
Section 6.3.3 of \cite{Ganb}.)

In any case, the Halphen or Schwarz differential equations can be used to compute arbitrarily
many terms of Fourier or Taylor expansions of automorphic forms (on the third author's 
homepage one can find computer code written in singular \cite{GPS01} and the first few coefficients of $t_1,t_2,t_3,J_\mft$ at $\I\infty$). From these expansions we are led to 
the conjectures (and results) gathered below. 

We will find a deep connection to the arithmeticity (or otherwise) of $\Gamma_\mft$,
and the integrality of those coefficients. This is hardly surprising. If a Fuchsian group
has at least 1 cusp (as we've been assuming), then the definition of arithmeticity
can be taken to be that it contains some conjugate of some congruence subgroup $\Gamma(N)$.
 By a theorem of Margulis \cite{Marg}, a Fuchsian group is arithmetic iff the commensurator
 $$\mathrm{comm}(\Gamma):=\{\gamma\in\mathrm{PSL}(2,\bbR)\,:\,\gamma\Gamma\gamma^{-1}\
 \mathrm{is\ commensurable\ with}\ \Gamma\}$$
is dense (recall that $\Gamma_1,\Gamma_2$ are \textit{commensurable} iff $\Gamma_1\cap
\Gamma_2$ has finite index in both $\Gamma_i$). More precisely, when $\Gamma$ is non-arithmetic, comm$(\Gamma)$ is itself
a Fuchsian group of the first kind, in fact the largest containing $\Gamma$. On the other hand,
if $\Gamma$ contains some $\Gamma(N)$ then any $\gamma\in\mathrm{GL}^+(
2,\bbQ)$ (or rather its projection to PSL$(2,\bbR)$) will lie in comm$(\Gamma)$. The
relevance of the commensurator is that $\gamma\in\mathrm{comm}(\Gamma)$ directly
yields Hecke operators for $\Gamma$. Given enough Hecke operators, the arithmeticity
of coefficients will follow.

It is easy to see directly that, for the non-arithmetic triangle groups, something goes wrong
with standard Hecke theory. Recall that the basis of Theorems 2 and 4 look like
$$
f(\tau)=\sum_{n=0}^\infty a_n\, q_3^n,\quad  a_n=r_n\,\alpha_3^n,
$$
where $r_n\in\bbQ$ and $q_3=e^{\frac{2\pi\I\tau}{h}}$.  \cite{wol83} proved that $\alpha_3$ is  transcendental,
but that implies that  $a_na_m\neq a_{mn}$ whenever $m,n>2$. Nor can we get multiplicativity
if we absorb the $\alpha_3$ into $q_3$. For weight $k$ cusp forms for any Fuchsian group,
we have the bound $a_n=O(n^{k/2})$ \cite{Miy}. But this means that the $r_n$ 
increase or decrease exponentially (depending on whether or not $|\alpha_3|<1$), which is again incompatible
with $r_nr_m=r_{nm}$ for sufficiently large $m,n$.

\subsection{Coefficients at the cusps}
Fix a hyperbolic type $\mft=(m_1,m_2,\infty)$. We do not require here that $m_1\le m_2$;
the case where $m_1$ or $m_2$ is infinite is included in the formulas below using the
aforementioned convention about the value of polynomials at $\infty$.
Consider first the Fourier coefficients $c_n
=c_{n;\mft}$ of  \eqref{jmnqexp}.
Note that the Euclidean types $(2,2,\infty)$ and (formally) $(1,\infty,\infty)$ correspond to
polynomial solutions $\widetilde{q}_3^{-1}+\frac{1}{2}+\frac{1}{16}\widetilde{q}_3$ and
$\widetilde{q}_3^{-1}$ respectively of \eqref{schw}. This means that $c_n$ vanishes when $m_1=m_2=2$ $\forall n\ge 2$,  and also $c_n$ vanishes at $m_1=1,m_2=\infty$ $\forall n\ge 0$, and hence
\begin{align} \label{24may2012}
&&c_n=\frac{(m_1^2-4)P_{1;n}(m_1^2,m_2^2)+(m_2^2-4)P_{2;n}(m_1^2,m_2^2)}{ (m_1^2m_2^2)^{n+1}\,Q_n}\,,\qquad n\geq 2\,,\\  \label{24may2012b}
&&c_n=\frac{(m_1^2-1)P'_{1;n}(m_1^2,m_2^2)+m_2^2P'_{2;n}(m_1^2,m_2^2)}{(m_1^2m_2^2)^{n+1}\,Q'_n}\,,\qquad n\geq 1,\
\end{align}
where $Q_n,Q'_n\in \N$ and $P_{i;n},P'_{i;n}$ are type-independent polynomials with integral coefficients and
total degree $\leq n-1$. \eqref{24may2012} generalizes to any type $(m_1,m_2,\infty)$ 
the observation of Akiyama \cite{aki92} described below, and \eqref{24may2012b} seems
completely new. Note that it would be reasonable to absorb $(m_1^2m_2^2)^{n}$ into ${\tilde q}_3$, at least when $m_1,m_2$ are both finite, and indeed this gives the $\hat{q}$ used in
Section 3.

A more interesting symmetry is that for $n\ge 1$, 
\begin{equation}c_{n;(m_1,m_2,\infty)}=(-1)^{n+1}c_{n;(m_2,m_1,\infty)}\,.\end{equation}
To prove this, first identify $\Gamma_{(m_2,m_1,\infty)}$ as a conjugate of $\Gamma_{(m_1,m_2,\infty)}$,
and then use this to express $J_{(m_2,m_1,\infty)}$ in terms of $J_{(m_1,m_2,\infty)}$.

Some of this had already been worked out for the Hecke groups $\Gamma_{(2,m,\infty)}$. In particular,  Lehner  
\cite{leh54} and especially  Raleigh \cite{ral62} worked from the Schwarz equation, obtaining
\eqref{a-3} in this special case as well as (\ref{24may2012})  without the $m^2-4$ factor.
For $n\geq 2$ and again only for the Hecke groups, Akiyama \cite{aki92} showed that $c_n$ is a polynomial divisible by 
$m^2-4$. He also showed that
the prime divisors of $Q_n$ are not greater that $n+1$. This follows immediately from the recursion given by the Halphen differential equation, where at the $n$-th step
of the recursion we divide by $n^2(n-1)$, see \S\ref{proofshalphen}. 
Leo in his PhD thesis \cite{leo08} proved 
that $c_n$ can be written as $\frac{C_n}{D_n(2^6m^2)^{n+1}}$, where $C_n, D_n\in \Z$ are coprime and 
$D_n$ has no prime factor 
of the form $p\equiv 1\pmod {4m}$.
He made also a precise conjecture about the prime factors of $D_n$. As with all these people,
he focussed exclusively on the Hecke groups $\Gamma_{(2,m,\infty)}$.

A major conjecture, now attributed to Atkin and Swinnerton-Dyer \cite{AS}, states that if $f$
is a modular form of weight $k\in \frac{1}{2}\bbZ$ for some subgroup $\Gamma$ of $\Gamma(1)$,
and the Fourier coefficients are algebraic integers, then $\Gamma$ (if it is chosen maximally)
contains a congruence subgroup. See e.g. \cite{LL} for a review.  Scholl  \cite{scholl} has proved
that when $\Gamma$ is a subgroup of  $\Gamma(1)$, there is an integer $N$ and a scalar
multiple $\widetilde{q}$ of $q=e^{2\pi\I\tau}$ such that 
the space of modular forms for $\Gamma$ of each weight $k\in \frac{1}{2}\bbZ$ has a basis with
$\widetilde{q}$-expansion coefficients which are algebraic integers when multiplied by
some power of $N$. We have  $N=1$ if (and conjecturally only if) $\Gamma$ contains a congruence subgroup, i.e. is arithmetic. In other words, we know that at most finitely many distinct primes can appear
in the denominators of modular forms for subgroups of $\Gamma(1)$. On the other hand,
 when $\Gamma$ is not commensurable with $\Gamma(1)$, one would expect 
infinitely many distinct primes in the denominators.

Our observations are compatible with these conjectures. Recall from Section 4 the 
 $9$ arithmetic triangle groups with at least one cusp: namely those of type
\begin{equation}
 \label{take}
(\infty,\infty,\infty), (2,3,\infty), (3,3,\infty),(m,\infty,\infty),
(2,2m,\infty),(2m,2m,\infty),\end{equation} 
for $m=2,3$. This also coincides with the list of all triangle groups conjugate
 to a group commensurable with $\Gamma(1)$. All 9 of those (up to conjugation)
 contain a congruence subgroup, as they must. In Section 4 we recovered the classical result that
  in these cases the algebra of modular  forms for $\Gamma_\mft$ is defined over $\Z$. By that
  we mean that there is a rescaling $Q$ of $q_3$, and some modular forms $f_i\in\bbZ[[Q]],\ i=1,2,\ldots$,  such 
  that the algebra of all holomorphic modular forms for $\Gamma_\mft$ is $\bbC[f_i,\ i=1,2,\ldots]$. 

The algebra of automorphic forms for the hyperbolic triangle group $\Gamma_{(m_1,m_2,\infty)}$ is defined over $\Z$ if and only if 
the triangle group is arithmetic. The only if part of this affirmation  is classical, and was reproved in \S 4. The other direction
has been recently proved by the last two authors. For the non-arithmetic case we are also
able to prove that infinitely many primes do not appear  in any denominators  of the coefficients of 
$t_i,\ i=1,2,3$ and $J_\mft$.  We are led to the following conjecture experimentally:
\begin{conj}
\label{15may2013}
For any non-arithmetic hyperbolic triangle group of type $(m_1,m_2,\infty)$,  infinitely many primes 
appear  in the denominators  of the coefficients of $t_i,\ i=1,2,3$ and $J_\mft$ 
at the infinite cusp.
\end{conj}
For non-arithmetic $\Gamma_\mft$ with $2\leq m_1\leq m_2\leq 30$ (and several other $m_i$
chosen randomly),  we looked
at all denominators for terms up to $q^{182}$. The distribution of primes which appear,
compared with those which do not, seem to be similar.
We also observe  that for each prime $p\not= 2$, $t_i(p \hat{q}),\ i=1,2,3$ has no $p$ in the denominators of its coefficients. 
This can be easily seen from  the recursion
given by the Halphen differential equation, see \S\ref{proofshalphen}. 
More precisely,
let  $p$ be a prime and $f$ be an automorphic form for $\Gamma_{(m_1,m_2,\infty)}$. 
Define  $m_{n,p}(f)$ to be 
the power of $p$ in the denominator of $a_n$, where $f=\sum a_n\hat{q}^n$. Our data suggests
the conjecture  $\lim_{n\to \infty} \frac{m_{n,p}}{n}=0$.

The main thing responsible for this non-integrality is the coefficient $Q_n$ in the
denominator of  \eqref{24may2012}.  
We suspect that each prime appears in the prime decomposition of some $Q_n$. 
The reason is that in the recursion for calculating the coefficients of $\tilde{q}^n$ we divide by $n^2(n-1)$. Although a priori a prime $p$ could appear at $n=p$, 
we observe that it appears first
at $n=p+1$. Note that this observation does not imply  Conjecture \ref{15may2013}, since
 the denominator and numerator of $c_n$  in (\ref{24may2012}) may have common factors. 

The much simpler case of Hecke groups  is extensively analyzed by Leo in \cite{leo08}.
For completeness we review his findings. Consider the triangle group of type $(2,m,\infty)$. Write
$$c_n=\frac{C_n}{D_n2^{6n+6}m^{2n+2}}$$ where
$C_n,D_n\in\bbZ$ and gcd($C_n,D_n)=1$. 
Leo \cite{leo08} conjectured that a prime $p$ divides some $D_n$ for $n\ge 1$, iff $p\ne 2$,
$p$ doesn't divide $m$, and $p\not\equiv \pm 1$ (mod $m$). Moreover, he conjectures that
the smallest $n$ for which such a prime $p$ divides $D_n$, is $n=p^k-1$ for some $k$.



%
%
%
%


\subsection{Integrality at elliptic fixed-points}

Again, we propose studying these expansions because every triangle group has 3 special
$\Gamma_\mft$-orbits,
most of which are elliptic fixed-points. As already mentioned,
\cite{RVZ} has found some of these coefficients to be interesting.

Consider first $\Gamma_{(2,3,\infty)}=\Gamma(1)$. Recall the expansion \eqref{jmnqexp}.
The coefficients at $\tau=\I$ 
are \begin{equation}
a_2 =\frac{23}{54}\,,\  a_3=\frac{6227}{58320}\,,\  a_4=\frac{3319}{174960}\,,\  a_5=\frac{263489}{97977600}\,,\ a_6=\frac{1693777}{5290790400}\,,\ldots\label{Jcoeff}\end{equation}
Not only are these nonintegral, but the denominator seems to be growing without bound!
But as we shall see shortly, there is a simple explanation for this.

The coefficients at elliptic fixed-points are more accessible than the coefficients at cusps. 
In particular, choose any point $z=x+\I y\in\bbH$ of order $m\ge 1$ and let
$f(\tau)=j_z(k;\tau)q_z^{k/m}\sum c_nq_z^n$ be a weight-$k$ automorphic form 
(recall \eqref{eq:laurexp1}). Note that $q_z$ is \textit{not} rescaled here, so
that series will have radius of convergence exactly 1 (provided $f$ is holomorphic).
Incidentally, Cauchy-Hadamard constrains the growth of these $c_n$:
$\mathrm{lim\,sup}_{n\rightarrow\infty}|c_n|^{1/n}=1$, so they grow roughly like
the usual (unscaled) Fourier coefficients.
 
These coefficients $c_n$ are then computed by  \cite{RVZ,BGHZ}
\begin{equation}c_n=\partial_k^{nm}f(z)\frac{(4\pi y)^{mn}}{(mn)!}\,,\label{ellcoeff}\end{equation}
where $z=x+\I y$ and $\partial_k^n=\partial_{k+n-2}\circ\cdots\circ\partial_{k+2}\circ\partial_k$, for
the nonholomorphic modular derivative $\partial_kf=\frac{1}{2\pi\I}\frac{\mathrm{df}}{\mathrm{d}\tau}-\frac{kf(z)}{4\pi y}\,.$
The $mn$ arises because $q_z=(\cdots)^m$ is a power.
Hence in this sense we can think of these $c_n$ as Taylor coefficients. 
The reason for the terrible  denominators in \eqref{Jcoeff} is the $n!$ in   \eqref{ellcoeff}.

The important quantities
should be the derivatives of $f$, in other words we should multiply the $a_n$ by $n!$
(and rescale $q_z$). We find for $\Gamma(1)$ at $z=\I$ that $a_n(mn)!m^n$      are positive integers, with
a single 3 in the denominators. The analogous calculation for the other elliptic fixed-point yields
only positive integers. We expect:

\begin{conj} Consider any arithmetic triangle group $\Gamma_{(m_1,m_2,\infty)}$
and  any elliptic fixed-point $z\in\uhp$. Then the sequence $(m_1n)!m_2^na_n$ are 
strictly positive algebraic numbers with bounded denominators.
There should exist a basis for the space of weight $k$ holomorphic automorphic
forms whose coefficients at $z$ are algebraic integers when rescaled in this way.
\end{conj}

For $\mft=(3,3,\infty)$, the denominator for $J_\mft$ is bounded by 8, while for $(4,4,\infty)$ and $(6,6,\infty)$
the denominators are all 1. For $(2,4,\infty)$, the adjusted
$a_n$ have denominators bounded by 2, while the adjusted $b_n$ have at most 3 in the 
denominators. For $(2,6,\infty)$, the adjusted
$a_n$ have at most a 3 in the denominator, while the adjusted $b_n$ is integral.
The larger the order of the fixed-point, the greater the chance for integers, because the multipliers
become so big. Note that for an arithmetic triangle group $(m_1,m_2,\infty)$ it
suffices to compute the values $\partial_k^{nm}f(z)$ for the generators $f$, as $\partial_k$
is a derivation. 

For non-arithmetic types, the situation is less clear. For example, for $\mft=(2,5,\infty)$,
the adjusted $a_n$ has 5's appearing in the denominators to arbitrarily high powers, and the
only other prime appearing in a denominator is 2, with power at most 3. In this case
$a_n(2n)!5^{2n}$ has bounded denominators. On the other hand the adjusted $b_n$ is
integral. For $(m_1,m_2)=(2,7),(2,8),(3,7)$, $a_n(m_1n)!m_2^{n}$ has unbounded denominator
but $a_n(m_1n)!m_2^{2n}$
and $b_n(m_2n)!m_1^{2n}$ both have bounded denominators. All of these were verified
up to $n=35$, but because of recursive formulas for these coefficients, it shouldn't be difficult
to prove this.



\section{Periods and automorphic functions}
The Gauss hypergeometric functions are periods up to some $\Gamma$-factors. This means that we can write them as integrals
of algebraic differential forms over topological cycles. Looking in this way we can generalize automorphic functions beyond
their classical context of Hermitian symmetric domains and action of groups, see for instance 
Section \ref{14cases}. 
In this section we explain this idea.  
\subsection{Periods and Halphen}
In \cite{ho07-1} the third author has used integrals of the form $\int \frac{x^idx}{(x-t_1)^a(x-t_2)^b(x-t_3)^c},$ 
in order to establish various properties of Halphen differential equations so that  generalizations, for instance for arbitrary 
number of $x-t_i$ factors  in the integrand, become realizable. We can 
view these integrals as periods in the following sense. 
We define a new variable $y$ and consider the family of algebraic curves $C: y=(x-t_1)^a(x-t_2)^b(x-t_3)^c$ for rational numbers $a,b,c$. 
In this way hypergeometric functions up to some $\Gamma$ factors can be written as periods $\int_{\delta}\omega$, where $\omega$ is a differential form on  
$C$ without
residues around its poles and $\delta\in H_1(C,\Z)$, see \cite{shtswo}. Now, 
one can use the algebraic geometry machinery in order to study the coefficients
of $q$-expansions of automorphic functions, see for instance \cite{ksv}, or the arithmetic of hypergeometric functions, see \cite{shtswo}. 
In the next subsection we describe a similar situation with Calabi-Yau periods.

\subsection{Hypergeometric Calabi-Yau equations}
\label{14cases}
Let $\widetilde{X}$ be a Calabi-Yau threefold, and $\cM$ its moduli space of complex
structures. The (complex) dimension of $\cM$ equals the Hodge number $h^{2,1}$.
We are interested here in $h^{2,1}=1$, in which we can, in the simplest cases, identify $\cM$ with $\bbC\bbP^1\setminus\{0,1,\infty\}$, where the large complex structure point corresponds to $z=0$, the
conifold point to $z=1$, and the Landau-Ginzburg point to $z=\infty$. The simplest
example is the mirror family of the generic quintic hypersurface in $\bbC\bbP^4$, which
can be parametrized by $x_1^5+x_2^5+x_3^5+x_4^5+x_5^5-5z^{-1/5}x_1x_2x_3x_4x_5=0$
for $z\in\cM$.

A holomorphic family $\varpi(z)$ of holomorphic 3-forms will satisfy the \textit{Picard-Fuchs
equation}. This implies, for any 3-cycle $\gamma\in H_3(\widetilde{X};\bbC)$, the period
$\int_\gamma\varpi(z)$ will satisfy a generalized hypergeometric equation of order $2h^{2,1}+2=4$,
also called the Picard-Fuchs equation. Periods provide a (redundant) parametrization of
$\cM$. See e.g. \cite{Ho-book} for a systematic treatment of periods, Picard-Fuchs, and
related concepts. 

There are  precisely 23 integral variations of Hodge structure which can come from such $\widetilde{X}$ with $h^{2,1}=1$, corresponding to 14 different
Picard-Fuchs equations \cite{DM}. For simplicity we have selected in Table 2 one representative
for each equation. The Picard-Fuchs equation satisfied by the periods is
\begin{equation}\label{PicFuc}\delta^4-z\prod_{i=1}^4(\delta+a_i)=0\,,\end{equation}
where we write $\delta=z\mathrm{d}/\mathrm{d}z$, $a_3=1-a_2$, and $a_4=1-a_1$.
 Periods are subject to monodromy 
as we circle the special points in $\cM$, and these can be worked out explicitly.

\medskip\centerline{\textit{Table 2} The Picard-Fuchs equation and monodromy data of one-parameter models}
$$\vbox{\tabskip=0pt\offinterlineskip
  \def\tablerule{\noalign{\hrule}}
  \halign to 4in{\strut#
  \tabskip=0em plus1em &    
    \hfil#&\vrule#&\hfil#&\vrule#&    
\hfil#&\vrule#&\hfil#
\tabskip=0pt\cr
& $(a_1,a_2)$&& $(n_1, n_2)$ &&\hfil type\hfil\cr 
\tablerule
& $(\frac{1}{5},\frac{2}{5})$&&$(-4,-5)$&&\hfil$(5,\infty,\infty)$\hfil\cr 
\tablerule
&$(\frac{1}{6},\frac{1}{3})$&& $(-3,-3)$&&\hfil$(6,\infty,\infty)$\hfil\cr
\tablerule
&$(\frac{1}{8},\frac{3}{8})$&&$(-3,-2)$&&\hfil$(8,\infty,\infty)$\hfil\cr
\tablerule
& $(\frac{1}{10},\frac{3}{10})$&& $(-2,-1)$&&\hfil$(10,\infty,\infty)$\hfil\cr
\tablerule
& $(\frac{1}{4},\frac{1}{3})$&& $(-4,-6)$&&\hfil$(12,\infty,\infty)$\hfil\cr
\tablerule
&$(\frac{1}{6},\frac{1}{4})$& &$(-2,-2)$&&\hfil$(12,\infty,\infty)$\hfil\cr
\tablerule
&$(\frac{1}{12},\frac{5}{12})$& &$(-3,-1)$&&\hfil$(12,\infty,\infty)$\hfil\cr
\tablerule
&$(\frac{1}{4},\frac{1}{2})$&&$(-5,-8)$&&\hfil$(\infty,\infty,\infty)$\hfil\cr
\tablerule
&$(\frac{1}{3},\frac{1}{2})$&&$(-6,-12)$&&\hfil$(\infty,\infty,\infty)$\hfil\cr
\tablerule
&$(\frac{1}{6},\frac{1}{2})$&&$(-4,-4)$&&\hfil$(\infty,\infty,\infty)$\hfil\cr
\tablerule
&$(\frac{1}{3},\frac{1}{3})$&&$(-5,-9)$&&\hfil$(\infty,\infty,\infty)$\hfil\cr
\tablerule
&$(\frac{1}{4},\frac{1}{4})$&&$(-3,-4)$&&\hfil$(\infty,\infty,\infty)$\hfil\cr
\tablerule
&$(\frac{1}{6},\frac{1}{6})$&&$(-1,-1)$&&\hfil$(\infty,\infty,\infty)$\hfil\cr
\tablerule
&$(\frac{1}{2},\frac{1}{2})$&&$(-7,-16)$&&\hfil$(\infty,\infty,\infty)$\hfil\cr
\noalign{\smallskip}}}$$

In particular, 
fix an integral basis $\gamma_1,\ldots,\gamma_4$ of $H_3(\widetilde{X};\bbZ)$. 
This is done in \cite{AGKLSY,GL01} using Meijer functions. Collect the
periods into a column vector $\Pi(z)=(\int_{\gamma_1}\varpi(z),\, \ldots\,, \int_{\gamma_4}\varpi(z))^t$.
Then $\Pi(z)$ is a fundamental solution of \eqref{PicFuc}. In terms of the Meijer basis,
the monodromy matrices are:
\begin{equation}M_0=\left(\begin{matrix}1&0&0&0\\ -1&1&0&0\\ 1&-1&1&0\\ 0&0&-1&1\end{matrix}
\right)\ ,\ M_\infty=\left(\begin{matrix}n_1&1-n_1&n_2&1-n_2\\ -1&1&0&0\\ 1&-1&1&0\\ 0&0&-1&1\end{matrix}
\right)\, ,\end{equation}
and $M_1=M_0^{-1}M_\infty^{-1}$,
using the parameters $n_i$ of Table 2, where $M_0$ is the monodromy picked up along a small
counterclockwise circle going around $z=0$, etc.

Of course, these monodromy matrices together define a representation of $\pi_1(\cM)\cong 
\Gamma_{(\infty,\infty,\infty)}$. In 7 of the models we can do better though. The orders of $M_0$ and $M_1$
will always be infinite, but those of $M_\infty$ can sometimes be finite. If we let $m$ be
the order of $M_\infty$,  then this representation of $\Gamma_{(\infty,\infty,\infty)}$ factors
through to a representation of $\Gamma_{(m,\infty,\infty)}$. This type $(m,\infty,\infty)$ is
collected in the final column of Table 2. What we lose in going to a less familiar triangle group,
we gain in getting a much tighter representation. Indeed, \cite{BrTh} show that for the first model
in Table 2, and a few others, the monodromy representation of $\Gamma_{(m,\infty,\infty)}$ is faithful;
by contrast, the kernel of the natural surjection $\Gamma_{(\infty,\infty,\infty)}\rightarrow
\Gamma_{(m,\infty,\infty)}$ is a free group of infinite rank for any $m<\infty$. 
It is a remarkable fact that 7 of the cases in Table 2 are of infinite index (see \cite{BrTh}) 
and 3 cases are of finite index, see \cite{si-ven}.

\subsection{Vector-valued automorphic forms}
A solution to a Fuchsian differential equation  over a compact surface, can be interpreted
as a \textit{vector-valued automorphic form} (vvaf) simply by lifting the surface minus singularities
($\bbC\bbP^1\setminus\{0,1,\infty\}$ here) to its universal cover $\bbH$. This isn't a completely
trivial statement ---  see \cite{BaGa} for the general argument --- but in the special case
of these models this will be made manifest shortly.

\medskip\noindent\textbf{Definition.} \textit{Let $k\in 2\bbZ$, $\Gamma$ be a Fuchsian group,
and $\rho$ a group homomorphism $\Gamma\rightarrow\mathrm{GL}(d,\bbC)$. A} vector-valued
automorphic form $\bbX(\tau)$ \textit{of weight $k$ on $\Gamma$ with multiplier $\rho$ is a
meromorphic map $\bbX:\bbH\rightarrow\bbC^d$, meromorphic also at the cusps, obeying the
functional equation}\begin{equation}
(c\tau+d)^{-k}\,\bbX\left(\frac{a\tau+b}{c\tau+d}\right)=\rho\left(\begin{matrix}a&b\\ c&d\end{matrix}\right)\,
\bbX(\tau)\,.\end{equation}

Choosing $\mft=(m,\infty,\infty)$ for either $m=\infty$ or any $m>0$ with $\gamma_3^m=1$,
$\bbX(\tau):=\Pi(J_\mft(\tau))$ is a vvaf of weight 0 for $\Gamma_\mft$, for multiplier which can
be identified with the monodromy of the Picard-Fuchs differential equation. This gives a 
modular interpretation for periods.

Let's be more explicit. Perhaps the simplest way to describe a vvaf $\bbX$ of weight $k$ and 
rank $n$ is to state a differential equation 
\begin{equation}D_k^n+f_{n-1}D_k^{n-1}+\cdots+f_0=0\label{vvmfDE}\end{equation}
satisfied by all components of $\bbX$, together with enough information to identify
which solution corresponds to each component. Here,  $f_j$ is an automorphic form for
$\Gamma_\mft$ of weight $2j$,
$D_k$ is the differential operator of Theorem \ref{lemma1}(ii), and 
$D_k^j=D_{k+2j-2}\circ \cdots\circ D_{k+2}\circ D_k$.

Recall the data for $(\infty,\infty,\infty)$ collected in Section 4.2.
We have $D_2\theta^4_2=(2\theta^4_3-\theta_2^8)/3$, $D_2\theta^4_3=(2\theta_3^4\theta_2^4-\theta_3^8)/3$, $D_0J_\mft=\theta_4^4
J_\mft$, $\Delta_\mft=\theta_2^4\theta_3^4\theta_4^4$.
Recall the parameters $a_1,a_2$ collected in Table 2.
The vvaf $\bbX(\tau)$ has rank 4 and weight 0 and corresponds to the differential equation 
\eqref{vvmfDE} with 
$$f_3=\frac{10B+8C}{3}\,,\ \ \ f_2=\frac{20B^2}{9}+BC(a_1^2+a_2^2-a_2-a_1+\frac{41}{9})+\frac{11C^2}{9}\,,$$
$$f_1= -\frac{20B^3}{27}+B^2C\frac{-2-2a_2-2a_1+2 a_1^2+2a_2^2}{3}+BC^2\frac{1+
12a_2^2+ 12a_1^2-12a_2-12a_1}{9}-\frac{C^3}{27}\,,$$ 
$$f_0=C^3B(a_1^2a_2-a_1^2a_2^2-a_1a_2+a_1a_2^2)$$
where we write $A=\theta_3^4,B=\theta_2^4,C=\theta_4^4=A-B$.
This looks more complicated because it is a uniform formula for all $a_i$.

The solutions all have an expansion $\sum_n c_n(\tau)q^{n/2}$  and each coefficient $c_n(\tau)$ is a polynomial of degree at most 3 in $\tau$.
We can identify which solution to call $\bbX_1,\bbX_2,\bbX_3,\bbX_4$ -- these form a basis of the solution space, and
the components of a vvmf of weight 0 for $\Gamma(2)$. We know everything about these vvmf, e.g. their
multiplier (i.e. to which representation of $\Gamma(2)$ they correspond), their local expansions at
each of the 3 cusps 0,1,$\infty$, etc. The components lie in $\bbQ[[\sqrt{q}]]$ but not $\bbZ[[\sqrt{q}]]$. 
So what we lose in the simplicity of the local expansions, we gain in the simplicity of the functional
equations (which just involve the usual M\"obius transformations defining $\Gamma(2)$).

The $\Gamma_{(m,\infty,\infty)}$ expressions should have some advantages, since that is \textit{really} the group
doing the acting --- $\Gamma(2)$ is a bit of a formal trick. We will provide those expressions 
elsewhere. But the uniformity and familiarity of $\Gamma(2)$
of course has its advantages too. 
This gives an answer to the question: what is a modular interpretation for the Calabi-Yau threefold periods? An alternate
answer to this question generalizes the algebraic geometric definition of (quasi-)modular forms and the relation of the 
Halphen differential equation with the Gauss-Manin connection to the families of Calabi-Yau varieties, see \cite{ho20}. The relation between
these two approaches  is discussed in the next subsection. In future work we will
reinterpret questions involving periods into the automorphic language and explore
whether this sheds any new light on them.

\subsection{Periods and modular-type functions}
The  most important modular object arising from the periods of Calabi-Yau varieties  is the Yukawa coupling. 
Let $\psi_0=1+O(z)$  and $\psi_1:=\psi_0\ln(z)+O(z)$ be respectively the holomorphic and 
logarithmic solutions of the hypergeometric equation (\ref{PicFuc}). The Yukawa coupling
$Y:=n_0\frac{\psi_0^4}{( \psi_0\theta\psi_1-\psi_1\theta\psi_0)^3(1-z)}$ is holomorphic at $z=0$ and so it can be written
in the Calabi-Yau mirror map $q=e^{\frac{\psi_1}{\psi_0}}$:
$$
Y:= n_0+\sum_{d=1}^\infty n_d d^3\frac{q^d}{1-q^d}
$$
Here, $n_0:=\int_M\omega^3 $, where $M$ is the  $A$-model Calabi-Yau threefold of mirror symmetry
and $\omega$ is the K\"ahler $2$-form of $M$ (the Picard-Fuchs equation (\ref{PicFuc}) is satisfied by the periods of the B-model Calabi-Yau threefold).
The numbers $n_d$ are supposed to count the number of rational curves of degree $d$ in a generic $M$. For the first
item in Table 2 the first  coefficients $n_d$ are given by $n_d=5, 2875, 609250, 317206375,\cdots$.

The field  generated by  $z,\delta^{i}\psi_0,\ i=0,1,2,3, \theta\psi_1-\psi_1\theta\psi_0,\psi_0\theta^2\psi_1-\psi_1\theta^2 \psi_0$ 
over $\C$ and written in the coordinate $q$, has
many common features with the field generated by quasi-automorphic forms for the group $\Gamma$ generated by $M_0$ and $M_\infty$, see \cite{ho20}. This includes functional equations with
respect to $\Gamma$, the corresponding Halphen equation and so on. However, note that the former field is of transcendental degree 
$3$, whereas this new field is of transcendental degree $7$.    
This gives
a second modular interpretation of the periods of Calabi-Yau varieties.

\section{Proofs}
This section contains the proof of the theorems announced earlier.

\subsection{Proof of Theorem \ref{first}}

Fix any hyperbolic $\mft=(m_1,m_2,\infty)\ne (\infty,\infty,\infty)$ (the extreme case
$(\infty,\infty,\infty)$ can be verified using case $\mathbf{\infty}^4$ in the appendix or by recalling 
familiar facts from the Fuchsian group $\Gamma(2)$). The hypergeometric
parameters $\tilde a,\tilde b,\tilde c$ are related to the angular ones $v_i=1/m_i$ via:
\begin{equation*}
\tilde a=\tilde b=(1-v_1-v_2)/2\,,\ \tilde c=1-v_1\,.
\end{equation*}

Let's begin with the derivation of the fundamental domain and generators of $\Gamma_\mft$.
Define the Schwarz  function
\begin{equation}\label{schw2}\phi(z)=\mu\frac{u_2(z)}{u_1(z)}=\mu
\frac{z^{1-\tilde c}\,F(\tilde a-\tilde c+1,\tilde b-\tilde c+1,2-\tilde c;z)}{F(\tilde a,\tilde b,\tilde c;z)}\,,\end{equation}
where $u_i$ are the independent solutions to the hypergeometric equation  given in \eqref{solns-hge} and the scale factor $\mu$ is \cite{Har}
\begin{equation}\label{mu}\mu=
\frac{\sin(\pi\,(\tilde c-\tilde a))}{\sin(\pi\,\tilde a)}\frac{\Gamma(\tilde a-\tilde{c}+1)^2\,\Gamma(\tilde{c})}{\Gamma(\tilde{a})^2\,\Gamma(2-\tilde{c})}\,,\end{equation}
and is chosen
to fit the target into the unit disc.
Then $\phi(z)$ maps the upper hemisphere of $\bbC\bbP^1\setminus\{0,1,\infty\}$ biholomorphically onto the (open) hyperbolic triangle in the Poincar\'e unit disc 
with vertices $\phi(0)=0, \phi(1)=\xi_2,\phi(\infty)=e^{\pi\I v_1}=-\zeta_1^{-1}$,  where 
${\xi}_2=\sin(\pi \tilde a)/\sin(\pi(\tilde c-\tilde a))$.
These values are calculated directly from the data in Appendix A. 
We can extend $\phi$ to all of $\bbC\bbP^1$ by reflecting in the real axis (so the triangle
becomes a quadrilateral), and we can
make $\phi$ into a multivalued function onto the full Poincar\'e disc by reflecting in the sides of that
quadrilateral.
The local expansion of $\phi$ about $z=0$ of course is obtained from \eqref{hyperg}, while those
about $z=1$ and $\infty$ are obtained from the formulas in cases $\mathbf{\infty}^0,\mathbf{\infty}^1,\mathbf{\infty}^2$ of  Appendix A.

We can map the unit disc to the upper half plane via
\begin{equation}\label{schw3}\tau(z)=\frac{\phi(z)+\zeta_1}{\zeta_1\phi(z)+1}\,.\end{equation}
It is easy to verify that $\tau(z)$ maps the unit disc to $\uhp$, and sends $z=0,1,\infty$ to
$\zeta_1,\zeta_2,\I\infty$. This means  the
normalized Hauptmodul $J_\mft(\tau)$ is related to the inverse map $z(\tau)$ by $J_\mft=1-z$. 
The monodromy of \eqref{gauss} directly yields the action $\left({\alpha\atop \gamma}{\beta\atop \delta}\right).\phi=(\delta\phi+\gamma)/(\beta\phi+\alpha)$, which up to conjugation reduces to
the action of $\Gamma_\mft$ on $\tau$.
The values of $\alpha_i$ (and $h_3$) can be computed from the $z=1,0,\infty$ asymptotics given in Appendix A, but were already computed in \cite{wol83}.
\eqref{schw} is simply the Schwarz equation \eqref{schwder} expressed in local coordinates.

\subsection{Proof of Theorem \ref{lemma1}}
\smallskip Now turn to Theorem \ref{lemma1}. Write $\mathfrak{m}_k$ for the space of
holomorphic automorphic forms of weight $k$.

The divisor div$\,f$ of a meromorphic automorphic form $f$ ($f$ not identically 0) is defined to be 
the formal (and finite) sum $\sum \mathrm{ord}_{[z]}(f) \,[z]$ where $[z]$ denotes the orbit 
$\Gamma_\mft z$. The degree of div$\,f$ for any such $f$ of weight $k$ for a triangle
group of type $(m_1,m_2,m_3)$ is (see Theorem 2.3.3 of \cite{Miy} for a generalization)
\begin{equation}\label{degree}
\mathrm{deg}(\mathrm{div}\, f)=\frac{k}{2}\left(1-\frac{1}{m_1}-\frac{1}{m_2}-\frac{1}{m_3}\right)
\,.\end{equation}
  By the classical argument, $\dot{J}_\mft$
is an automorphic form for $\Gamma_\mft$ of weight 2, since $J_\mft$ has weight 0.  Clearly, the 
only poles of $\dot{J}_\mft$ are at the points  in $[\I\infty]$, where we have a simple pole.
Also, $\dot{J}_\mft$ has zeros at any other cusp (with order $\ge 1$) and at
elliptic points $\zeta_{i}$ (with order $\ge 1-1/m_i$). That these
orders are equalities, and that $\dot{J}_\mft$ has no other zeros, follows from the formula
for the degree of the divisor.

It is manifest from the formula for $f_k$ that   is an automorphic form of weight $k$,  holomorphic
everywhere in $\uhp_\mft$ except possibly at $[\I\infty]$.
 Note that for automorphic forms $f,g$ of fixed weight, the orders of $f$ and $g$ at any point will differ
  by an integer, and thus 
 the order of $f_k$ at each point $\not\in[\I\infty]$
is the minimum possible for $f\in\frak{m}_k$.

The quantity $d_k$ equals the order of $f_k$ at $\I\infty$. 
If $d_k\ge 0$ then for each $0\le l\le d_k$, $f_kJ_\mft^l$ is holomorphic at $\I\infty$ (hence lies in
$\frak{m}_k$). In this case, for any $g\in \frak{m}_k$,
$g/f_k$ will be an automorphic function holomorphic everywhere in $\uhp_\mft$ except possibly
at $\I\infty$. This means $g/f_k$ must equal
some polynomial in $J_\mft$ of degree $\le d_k$. Thus the $f_k J_\mft^l$ span 
$ \frak{m}_k$.
On the other hand, if $d_k<0$, then $\frak{m}_k=0$ (again, look at $g/f_k$ for any
$g\in\frak{m}_k$).

Consider now the generators of the algebra of holomorphic modular forms. 
Type $(\infty,\infty,\infty)$
can be obtained by recalling what is known for $\Gamma(2)$. Suppose first that
 $m_1<m_2=\infty$. Choose any $k\ge 0$ and write $k=k'+l m_1$ for $0\le k'<m_1$ and $l\in\bbZ$.
Note that $f_{2k'}\,f_{2m_1}^j$ has weight $2k$ and has order $1-k'/m_1$ (the smallest
possible in $\mathfrak{m}_{2k}$) at $\zeta_1$.
Then, given $f\in \mathfrak{m}_{2k}$, a constant $c$ can be found so that $f-c\,f_{2k'}\,f_{2m_1}^j$
will have order $\ge1$ at $\zeta_1$. Since $f_2$ has order $1-1/m_1,0,0$ at $\zeta_1,
\zeta_2,\zeta_3$ respectively, $(f-c\,f_{2k'}\,f_{2m_1}^j)/f_2\in\mathfrak{m}_{2k-2}$.
Thus by induction, $f_2,\ldots,f_{2m_1}$ generate all of $\mathfrak{m}_{2k}$, for any $k$.

The proof for  $m_2<\infty$ is similar. Define $f^{(1)}_{2l}:=f_{2l}J_\mft^{d_{2l}}$ (minimal 
possible order at $\zeta_1$ and $\I\infty$,
maximal at $\zeta_2$, in $\mathfrak{m}_{2l}$). Choose any $f\in\frak{m}_{2k}$ for $k\ge 0$, and 
write  $k=k_i+l_im_i$ for $i=1,2$ where $0\le k_i<m_i$ and $l_i\in\bbZ$. Then it is possible to 
find constants $c_i$ so that
$g:=f-c_2\,(f_{2m_2})^{l_2}\,f_{2k_2}-c_1\,(f_{2m_1})^{l_1}\,f_{2k_1}^{(1)}$ has order $\ge 1$ at both $\zeta_1,\zeta_2$.
This means $g/f_4\in\frak{m}_{k-4}$, so the result follows by induction on $k$.

As defined, $\Delta_{\mft}$ is manifestly a weight $2L$ automorphic form with no zeros or poles anywhere
except possibly at $[\I\infty]$. In fact, since $J_{\mft}$ is a Hauptmodul, the
order of $J_\mft(\tau)-J_\mft(\zeta_i)$ at $\zeta_i$ equals 1,  which gives us the formula for
$n_\Delta$. That value
is proportional to the area of a fundamental domain of $\Gamma_\mft$ (see e.g. \cite{Miy}), and so is strictly 
positive. Hence $\Delta_{\mft}$ vanishes at $\I\infty$.

The statement about holomorphicity of $E_{2;\mft}$ is immediate from the properties of
$\Delta_{\mft}$. The functional
equation for $E_{2;\mft}$ follows directly from that of $\Delta_{\mft}$, and
the vanishing of $E_{2;\mft}(\zeta_{j})$ at cusps $\zeta_{j}$ is a consequence of $\Delta_\mft(\zeta_{j})$ being
finite and nonzero there.

\subsection{Proof of Theorem \ref{15may2013-1}}
\label{proofshalphen}

The only new part of Theorem \ref{9feb2012} is (ii). Write $h=h_3$. Many
of their properties can be easily determined from those of the hypergeometric functions
collected in Appendix A. In particular, 
they are meromorphic functions in $\uhp$ with possible poles only at the $\Gamma_\mft$-orbits
of $\zeta_2$ and $\zeta_1$. Now, each $t_i$ is a function of
$\qh$, because $J_\mft$ is. 
Write $t_i=\sum_{n=0}^\infty t_{i,n}\qh^n$.
We see directly from (i) that, in vector form, 
$$
[t_{1,0},t_{2,0},t_{3,0}]=[0,-2\pi\I/h,0]\,
$$
(these are normalized differently in Theorem 3).
Comparing $\qh^n$ coefficients, for $n\geq 1$, we get a recursion:  
\begin{equation}
 \label{31maynight}
(M-nI_{3\times 3})[t_{1,n},t_{2,n},t_{3,n}]^\tr\in\bbQ^3[a,b,c][t_{i,m}]_{1\le i\le 3,0\le m<n}\,,
\end{equation}
where $\tr$ denotes transpose and
$$
M:=
\left(
\begin{array}{*{3}{c}}
\frac{m_{1}m_2+m_{2}-m_1}{2m_{1}m_2} & 0 & -\frac{m_{1}m_2+m_{2}-m_1}{2m_{1}m_2} \\
\frac{m_{1}m_2+m_{1}+m_2}{2m_{1}m_2} & 0 & \frac{m_{1}m_2+m_{1}+m_2}{2m_{1}m_2} \\
-\frac{m_{1}m_2-m_{2}+m_1}{2m_{1}m_2} & 0 & \frac{m_{1}m_2-m_{2}+m_1}{2m_{1}m_2}
\end{array}
\right).
$$
Note that 
$$
(M-I_{3\times 3})[t_{1,1},t_{2,1},t_{3,1}]^\tr=0\,,
$$
and so up to a constant $\nu'$,
{\small
$$
[t_{1,1},t_{2,1},t_{3,1}]=
\nu'\left[
\begin{array}{*{3}{c}}
-m_{1}^{2}m_2^{2}-m_{2}^{2}m_1+m_{2}m_1^{2}, & 
-m_{2}^{2}m_1-m_{2}^{2}+m_{2}m_1^{2}+m_1^{2}, & 
m_{1}^{2}m_2^{2}-m_{2}^{2}m_1+m_{2}m_1^{2}
\end{array}
\right]
$$
}
when $m_2<\infty$, while
$$[t_{1,1},t_{2,1},t_{3,1}]=
\left\{\begin{matrix}\nu'\left[
\begin{array}{*{3}{c}}
-m_1^{2}-m_1, & -m_1-1, & m_1^{2}-m_1
\end{array}
\right]&\mathrm{if}\ m_1<\infty=m_2\\
\nu'\left[
\begin{array}{*{3}{c}}
-1, & 0, & 1
\end{array}
\right]&\mathrm{if}\ m_1=m_2=\infty\end{matrix}\right.\,.
$$
(The rule is that the value of a polynomial $P(x)$ for $x=\infty$ is the coefficient of the  monomial $x^n$ of highest degree in $P(x)$.)
We chose the constant $\nu'$ here so that these expressions are polynomial in $m_2$ and $m_1$. 
That $\nu'=\nu$ follows by computing the leading term of $t_1$.

Note that $\det(M-nI_{3\times 3})=-n^2(n-1)$ so the $n$th coefficients of  $t_i$ are well-defined
polynomials in $m_j$ for $n>1$. 
The factor of $2\pi\I/h$ and power of $\nu$ in \eqref{ti} follows from easy inductions. In order to see (\ref{31may2012}), we write \eqref{tatlis}  in the variables $x_1=(m_1m_2)^{-2}(t_1-t_3),\ x_2=\kappa_2^{-1}(t_2+1)$ and $x_3=\kappa_3^{-1}t_3$, and we get a similar recursion as in (\ref{31maynight}) with different $M$ 
such that $\det(M-nI_{3\times 3})\ne 0$  even for $m_2=0$ and $m_1=0$. 
\subsection{Proof of Theorem \ref{29feb2012}}
That the $t_i$ obey \eqref{aruquia} is clear from Theorem \ref{9feb2012} and the
automorphy of $J_\mft$.
We obtain from Appendix A that, when $m_2\ne \infty$, the zero and pole orders of the $t_i$'s at the $\zeta_j$'s
are given in the table below:  
\begin{center}
\label{7mar2012}
\begin{tabular}{ l | l l l}
            & $\zeta_1$& $\zeta_2$ & $\zeta_3$ \\ \hline
 $t_2-t_1$&  $m_1-1$ & $-1$ &$0$ \\
$t_3-t_2$ &$-1$& $m_2-1$ &$0$   \\
  $t_1-t_3$ &$-1$& $-1$ & $1$ \\
$t_1$&$-1$&$-1$&$1$\\
$t_2$&$-1$&$-1$&$0$\\
$t_3$       &$-1$& $-1$ &  $1$ 
\end{tabular}
\end{center}
while if $m_1<\infty=m_2$, the table becomes
\begin{center}
\begin{tabular}{ l |l l l}
            & $\zeta_1$&$\zeta_2$ & $\zeta_3$ \\ \hline
  $t_2-t_1$ & $m_1-1$&$0$ &$0$  \\
$t_3-t_2$ & $-1$&$1$ &  $0$ \\
  $t_1-t_3$ & $-1$&$0$ & $1$ \\
$t_1$       & $-1$&$0$ &  $1$\\
$t_2$       & $-1$&$1$ &  $0$\\
$t_3$       & $-1$&$1$ &  $1$ 
\end{tabular}
\end{center}
(Note however that the  orders of zeros for quasi-automorphic forms like $t_i$ are not  constant along orbits.)
This table makes it easy to verify the automorphic form identities given in Theorem \ref{29feb2012}(iii). For the identity involving $E_{2;\mathfrak t},t_1,t_2,t_3$ 
we must further calculate the residues of $t_i$'s at elliptic points $\zeta_i$'s. Theorem \ref{29feb2012}(iv) follows by similar pole order arguments as in the proof of  
Theorem \ref{lemma1} and the above tables.

 \appendix
%
%
%
%
 
\section{Hypergeometric functions: Basic formulas}
\label{3apr2012}

In this appendix we review some classical facts about the \textit{Gauss hypergeometric function}
(or \textit{series})\begin{equation}\label{hyperg}
F(\tilde a,\tilde b,\tilde c;z)={}_2F_1(\tilde a,\tilde b,\tilde c;z)=\sum_{n=0}^{\infty} \frac{(\tilde a)_n(\tilde b)_n}{(\tilde c)_n n!}z^n, \ \tilde c\not\in
\{0,-1,-2,-3,\ldots \}\,,\end{equation}
where  $(x)_n:=x(x+1)(x+2)\cdots (x+n-1)$, and its differential equation
\begin{equation}
\label{gauss2}
z(1-z)y''+(\tilde c-(\tilde a+\tilde b+1)z)y'-\tilde a\tilde by=0,
\end{equation}
which is called the \textit{hypergeometric} or \textit{Gauss equation}. A very complete reference is \cite{BaEr}, though it has typos.
In the following and throughout this paper, $\Gamma(z)$ denotes the gamma function
and the digamma $\psi(z)$ denotes its logarithmic derivative. The values of $\psi$
at rational $z$ (the only ones we need) were calculated by Gauss to be:
\begin{equation}\label{psirat}
\psi(m/n)=-\gamma-\ln\,n-\frac{\pi}{2}\cot(\pi m/n)+\sum_{k=1}^{n/2}{}'\cos(2\pi mk/n)\ln(2-2\cos(2\pi k/m))
\end{equation}
where $\gamma$ is Euler's constant and the prime means that for even $n$ the last term (namely, 
$k=n/2$) should be divided by 2. Another identity is useful:
\begin{equation}\label{psi2}
\psi(1-x)=\psi(x)+\pi\cot \pi x\,.\end{equation}

 The values $\tilde a,\tilde b,\tilde c$ of interest here are given at the beginning of Section 7.1
 and (more generally) Appendix B. 
As long as $\tilde c\not\in\bbZ$  (i.e. except for case $\mathbf{\infty}^3$ below), the solution 
space to  \eqref{gauss2} is spanned by
\begin{equation}u_1(z)=F(\tilde a,\tilde b,\tilde c;z)\,,\ \ u_2(z)=
z^{1-\tilde c}F(\tilde a-\tilde c+1,\tilde b-\tilde c+1,2-\tilde c;z)\,.\label{solns-hge}\end{equation}
We need to understand what $u_i(z)$ looks like about $z=1$ and $z=\infty$, in order to
understand the local expansions of the automorphic forms of $\Gamma_{(m_1,m_2,m_3)}$ about
all cusps and elliptic fixed-points. Closely related to this, we need to understand the monodromy 
of \eqref{gauss2}
in order to explicitly identify the automorphic forms associated to $\Gamma_{(m_1,m_2,m_3)}$
(it is easy to identify them up to a conjugate of $\Gamma_{(m_1,m_2,m_3)}$, but we want to
pin down that conjugate). These formulas only depend on the number of cusps, i.e. the
number of $m_i$ which equal $\infty$. We will require here (without loss of generality)
that $m_1\le m_2\le m_3\le\infty$.

\medskip\noindent\textbf{Case $\mathbf{\infty}^0$:} \textit{No cusps, i.e. $m_3<\infty$.}\smallskip

This corresponds to all of $\tilde a,\tilde b,\tilde c,\tilde c-\tilde a-\tilde b,\tilde a-\tilde b$ 
being nonintegral.
Analytic continuation for $|\arg(1-z)|<\pi$ resp. $ |\arg(-z)|<\pi$ is: \begin{eqnarray*}
&&u_1(z)= \frac{\Gamma(\tilde c)\Gamma(\tilde c-\tilde a-\tilde b)}{\Gamma(\tilde c-\tilde a)\Gamma(\tilde c-\tilde b)}F(\tilde a,\tilde b,\tilde a+\tilde b-\tilde c+1;1- z)\\&&\qquad\qquad+
\frac{\Gamma(\tilde c)\Gamma(\tilde a+\tilde b-\tilde c)}{\Gamma(\tilde a)\Gamma(\tilde b)}   (1-z)^{\tilde c-\tilde a-\tilde b}F(\tilde c-\tilde a,\tilde c-\tilde b,\tilde c-\tilde a-\tilde b+1;1-z)  \\
&&\ \ \qquad= \frac{\Gamma(\tilde c)\Gamma(\tilde b-\tilde a)}{\Gamma(\tilde b)\Gamma(\tilde c-\tilde a)} (-z)^{-\tilde a}F(\tilde a,1-\tilde c+\tilde a,1-\tilde b+\tilde a;z^{-1})\\&&\qquad\qquad+ 
\frac{\Gamma(\tilde c)\Gamma(\tilde a-\tilde b)}{\Gamma(\tilde a)\Gamma(\tilde c-\tilde b)}   (-z)^{-\tilde b}F(\tilde b,1-\tilde c+\tilde b,1-\tilde a+\tilde b;z^{-1})\,,\\
&&u_2(z)= \frac{\Gamma(2-\tilde c)\Gamma(\tilde c-\tilde a-\tilde b)}{\Gamma(1-\tilde a)\Gamma(1-\tilde b)}F(\tilde a,\tilde b,\tilde a+\tilde b-\tilde c+1;1- z)\\&&\qquad\qquad+
\frac{\Gamma(2-\tilde c)\Gamma(\tilde a+\tilde b-\tilde c)}{\Gamma(\tilde a-\tilde c+1)\Gamma(\tilde b-\tilde c+1)}   (1-z)^{\tilde c-\tilde a-\tilde b}F(\tilde c-\tilde a,\tilde c-\tilde b,\tilde c-\tilde a-\tilde b+1;1-z)  \\
&&\ \ \qquad= -e^{-\pi\I\tilde c}\frac{\Gamma(2-\tilde c)\Gamma(\tilde b-\tilde a)}{\Gamma(\tilde b-\tilde c+1)\Gamma(1-\tilde a)} (-z)^{-\tilde a}F(\tilde a-\tilde c+1,\tilde a,1-\tilde b+\tilde a;z^{-1})\\&&\qquad\qquad-e^{-\pi\I\tilde c}
\frac{\Gamma(2-\tilde c)\Gamma(\tilde a-\tilde b)}{\Gamma(\tilde a-\tilde c+1)\Gamma(1-\tilde b)}   (-z)^{-\tilde b}F(\tilde b-\tilde c+1,\tilde b,1-\tilde a+\tilde b;z^{-1})\,.
\end{eqnarray*}
From this we obtain
the monodromy matrices (in terms of the basis $u_1,u_2$) 
for small counterclockwise circles about $z=0$, $z=1$, $z=\infty$:
\begin{align}
M_0=&\,\left(\begin{matrix}1&0\\ 0&e^{-2\pi\I \tilde{c}}\end{matrix}\right)\ ,&\\
M_1=&\,\left(\begin{matrix}\frac{\xi s(\tilde a)\,s(\tilde b)-s(\tilde c-\tilde a)\,s(\tilde c-\tilde b)}{s(
\tilde c)\,s(\tilde c-\tilde a-\tilde b)}&\frac{\pi\,(\xi-1)\Gamma(1-\tilde c)
\Gamma(2-\tilde c)}{s(\tilde c-\tilde a-\tilde b)\Gamma(1-\tilde a)\Gamma(1-\tilde b)\Gamma(\tilde a-\tilde c+1)\Gamma(\tilde b-\tilde c+1)}\\
\frac{\pi\,(\xi-1)\Gamma(\tilde c-1)\,\Gamma(\tilde c)}{s(\tilde c-\tilde a-\tilde b)\,\Gamma(\tilde c-\tilde a)\,\Gamma(\tilde c-\tilde b)\,\Gamma(\tilde a)\,\Gamma(\tilde b)}&\frac{s(\tilde a)\,s(\tilde b)-\xi s(\tilde c-\tilde a)\,
s(\tilde c-\tilde b)}{s(\tilde c)\,s(\tilde c-\tilde a-\tilde b)}
\end{matrix}\right)\,,&\end{align}
and  $M_\infty=M_1^{-1}M_0^{-1}$, where here $\xi=e^{\pi\I(\tilde c-\tilde a-\tilde b)}$ and $s(x)=\sin(\pi x)$.

\medskip\noindent\textbf{Case $\mathbf{\infty}^1$:} \textit{Exactly one cusp, i.e. $m_2<m_3=\infty$.}\smallskip

This means $\tilde a=\tilde b$, and all of $\tilde a,\tilde c,\tilde c-2\tilde a$ are nonintegral.
Analytic continuation to $z=1$ is as in case $\mathbf{\infty}^0$, but to $z=\infty$ is given by:
\begin{eqnarray*}
&u_1(z)=&\frac{(-z)^{-\tilde a}\Gamma(\tilde c)}{\Gamma(\tilde a)\Gamma(\tilde c-\tilde a)}
\sum_{n=0}^\infty\frac{(\tilde a)_n(1-\tilde c+\tilde a)_n}{n!n!}(\ln(-z)+2\psi(1+n)-\psi(\tilde a+n)-\psi(\tilde c-\tilde a-n))z^{-n}\\
&u_2(z)=&\frac{-e^{-\pi\I\tilde{c}}(-z)^{-\tilde a}\Gamma(2-\tilde c)}{\Gamma(1-\tilde a)\Gamma(\tilde a-\tilde c+1)}
\sum_{n=0}^\infty\frac{(\tilde a)_n(1-\tilde c+\tilde a)_n}{n!n!}(\ln(-z)+2\psi(1+n)\\
&&\qquad\qquad\qquad\qquad\qquad\qquad\qquad\qquad-\psi(\tilde a-\tilde c+n+1)-\psi(1-\tilde a-n))z^{-n}\,.
\end{eqnarray*}
Monodromy is given by the same matrices as in case $\mathbf{\infty}^0$.

\medskip\noindent\textbf{Case $\mathbf{\infty}^2$:} \textit{Exactly two cusps, i.e. $m_1<m_2=m_3=\infty$.}\smallskip

This means $\tilde a=\tilde b$ and $\tilde c=2\tilde a$, and both of $\tilde c,\tilde a$ are nonintegral.
Analytic continuation to $z=\infty$ is as in case $\mathbf{\infty}^1$, but  to $z=1$ is given by\begin{eqnarray*}
&&u_1(z)=\frac{\Gamma(2\tilde a)}{\Gamma(\tilde a)\Gamma(\tilde c)}\sum_{n=0}^\infty\frac{(\tilde a)_n(\tilde a)_n}{n!n!}(
2\psi(n+1)-2\psi(\tilde a+n)-\ln(1-z))\,(1-z)^n\\
&&u_2(z)=\frac{z^{1-2\tilde a}\Gamma(2-2\tilde a)}{\Gamma(1-\tilde a)\Gamma(1-\tilde a)}\sum_{n=0}^\infty
\frac{(1-\tilde a)_n(1-\tilde a)_n}{n!n!}(2\psi(n+1)-2\psi(1-\tilde a)\\ &&\qquad\qquad\qquad\qquad\qquad\qquad\qquad\qquad-\ln(1-z))\,(1-z)^n\,.
\end{eqnarray*}
Monodromy is again given by the same matrices as in case $\mathbf{\infty}^0$.

\medskip\noindent\textbf{Case $\mathbf{\infty}^3$:} \textit{Three cusps, i.e. $m_1=m_2=m_3=\infty$.}\smallskip

Then $\tilde a=\tilde b=1/2,\tilde c=1$. Take $u_1(z)=F(1/2,1/2,1;z)$ and
 \begin{equation*}
 u_2(z)=\I F(1/2,1/2,1;1-z)=\frac{\I}{\pi}\sum_{n=0}^\infty\frac{(1/2)_n(1/2)_n}{n!n!}(2\psi(1+n)-2\psi(1/2+n)-\ln(z))\,z^n\,,\end{equation*}
where the second equality is valid for $-\pi<\arg(z)<\pi$.
Analytic continuation of $u_1$ is as in case $\mathbf{\infty}^2$, but for $u_2$ is given by 
\begin{align*}u_2(z)&\,=\I F(1/2,1/2,1;1-z)&\,\\&\,=\frac{\I}{\pi} z^{-1/2}\sum_{n=0}^\infty\frac{(1/2)_n(1/2)_n}{n!n!}(2\psi(1+n)-2\psi(1/2+n)-\ln(z^{-1}))\,z^{-n}\,.&\end{align*}
The monodromy is
\begin{equation}M_0=\left(\begin{matrix}1&2\\ 0&1\end{matrix}\right)\ ,\ M_1=\left(\begin{matrix}1&0\\ -2&1\end{matrix}\right)\ ,\ M_\infty=\left(\begin{matrix}1&-2\\ 2&-3\end{matrix}\right)\ .\end{equation}


\section{Triangular groups without cusps}
\label{17may2013}
In this paper (and the applications we have in mind), we are interested in triangle groups with cusps, but the same calculations
work (though are messier) when there are no cusps, i.e. when all $m_i$ are finite. In this appendix
 we sketch the changes.  

 The equation \eqref{schw} becomes 
 \begin{equation}\label{schw22}-2{\dddot{J}_{\mft}\,\dot{J}_\mft}+{3}{\ddot{J}_{\mft}^2}
-n_z^{-2}\dot{J}_\mft^2=\dot{J}_{\mft}^4\left({1-v_2^2\over J_{\mft}^2}+{1-v_1^2\over (J_{\mft}-1)^2}+{v_1^2+v_2^2-v_3^2-1
\over J_{\mft}(J_{\mft}-1)}\right)\,.\end{equation}
For example,
 \begin{align*} c_0&=\frac{-1 +  \gamma_- + v_3^2}{2 (v_3^2 - 1)} \,,\qquad 
 c_1=\frac{(5 - 2 \gamma_+ - 3 \gamma_-^2 ) + (-6 + 2 \gamma_+) v_3^2 + v_3^4}{16(v_3^2-1)(v_3^2-4)}\,,\\
c_2&=\frac{(- 2 \gamma_-+\gamma_+ \gamma_-  + \gamma_-^3 ) + (2 \gamma_- - \gamma_+ \gamma_-) v_3^2}{6(v_3^2-9)(v_3^2-1)^2}\,,\\
c_3&=\frac{-31 + 76 \gamma_++ 690 \gamma_-^2  - 28 \gamma_+^2  - 404 \gamma_-^2  \gamma_+- 303 \gamma_-^4    }{128(v_3^2-16)(v_3^2-4)^2(v_3^2-1)^3}\\ & \quad+  \frac{100  - 244 \gamma_++ 88 \gamma_+^2   - 1052 \gamma_-^2+ 660 \gamma_-^2  \gamma_+  + 192 \gamma_-^4 }{128(v_3^2-16)(v_3^2-4)^2(v_3^2-1)^3}v_3^2\\ & \quad+\frac{-114  + 276 \gamma_+- 96 \gamma_+^2  + 390 \gamma_-^2 - 288 \gamma_-^2  \gamma_+  - 24 \gamma_-^4 }{128(v_3^2-16)(v_3^2-4)^2(v_3^2-1)^3}v_3^4\\ &\quad   +\frac{(52  - 124 \gamma_++ 40 \gamma_+^2  - 24 \gamma_-^2 + 32 \gamma_-^2  \gamma_+ ) v_3^6  +  (-7 + 16 \gamma_+ - 4 \gamma_+^2  - 4 \gamma_-^2 ) v_3^8}{128(v_3^2-16)(v_3^2-4)^2(v_3^2-1)^3}\,,\end{align*}
where  $\gamma_\pm=v_1^2\pm v_2^2$. 

The table in
Section 7.4, listing the orders of zeros and poles  for the solutions
of the Halphen system, generalizes to:
\begin{center}
\label{7mar2012b}
\begin{tabular}{ l | l l l}
            & $\zeta_3$& $\zeta_2$ & $\zeta_1$ \\ \hline
  $t_2-t_1$ & $-1$&$-1$ & $m_1-1$ \\
$t_3-t_2$ &$-1$& $m_2-1$ &  $-1$ \\
  $t_1-t_3$ &$m_3-1$& $-1$ & $-1$ \\
$t_1$&$-1$&$-1$&$-1$\\
$t_2$&$-1$&$-1$&$-1$\\
$t_3$       &$-1$& $-1$ &  $-1$ 
\end{tabular}
\end{center}
As before, a basis for the ring of automorphic forms consists of the monomials of the form
$$
(t_1-t_2)^p(t_2-t_3)^q(t_3-t_1)^r
$$
and the pole condition on the vertices implies that $p,q,r\geq 1$. 
The ring of holomorphic automorphic forms for the hyperbolic triangle group 
$\Gamma_{(m_1,m_2,m_3)}$ with the condition $m_1 \leq m_2\leq m_3<\infty$ is 
generated by holomorphic functions
\begin{align*}
 E_{p,q;\mft}^{(3)}&=(t_1-t_2)^p(t_2-t_3)^q(t_3-t_1), \quad k=p+q\,,\\
 E_{q,r;\mft}^{(1)}&=(t_1-t_2)(t_2-t_3)^q(t_3-t_1)^r, \quad k=q+r\,,\\
 E_{p,r;\mft}^{(2)}&=(t_1-t_2)^p(t_2-t_3)(t_3-t_1)^r, \quad  k=p+r\,.
\end{align*}
This list of generators is  finite because for example holomorphicity at $\zeta_3$ for $E_{p,q;\mft}^{(3}$ implies that $p+q\leq m_3-1$ and similarly for the rest.
The space of automorphic forms of weight $2k$ is of dimension 
 $k+1-\lceil \frac{k}{m_1}\rceil-\lceil \frac{k}{m_2}\rceil-\lceil \frac{k}{m_3}\rceil$.




  \newcommand\bibx[4]   {\bibitem{#1} {#2:} {\textit{#3}} {\rm #4}}

\end{document}